\documentclass[twoside,12pt]{amsart}
\usepackage[T1]{fontenc}
\usepackage{amsfonts}
\usepackage{amssymb}
\usepackage{amsthm}
\usepackage{amsmath}
\usepackage{graphicx}
\usepackage{url}
\usepackage[all]{xy}
\usepackage{lscape}
\usepackage{url}
\theoremstyle{plain}

\newcommand{\soc}{\mathrm{soc}}
\newcommand{\PSL}{\mathrm{PSL}}
\newcommand{\PSU}{\mathrm{PSU}}
\newcommand{\PSp}{\mathrm{PSp}}
\newcommand{\PGL}{\mathrm{PGL}}

\newcommand{\SO}{\mbox{\rm SO}}
\newcommand{\GL}{\mbox{\rm GL}}

\newcommand{\Aut}{\mbox{\rm Aut}}

\newcommand{\Aff}{\mbox{\rm Aff}}

\newcommand{\rcc}{\text{\rm RCC\ }}

\newenvironment{prf}{{\bf Proof.}}{\hfill $\Box$ \\[-1.0ex]}

\newtheorem*{theorem*}{Theorem}
\newtheorem{num}{Notation}[section]

\newtheorem*{define*}{Definition}

\newtheorem{thm}[num]{Theorem}
\newtheorem*{thm*}{Theorem}
\newtheorem{lem}[num]{Lemma}
\newtheorem*{lem*}{Lemma}
\newtheorem{prp}[num]{Proposition}
\newtheorem*{prp*}{Proposition}
\newtheorem{cor}[num]{Corollary}
\newtheorem*{cor*}{Corollary}

\newtheorem*{conj*}{Conjecture}

\newtheorem*{xmpl*}{Example}
\newtheorem{rem}[num]{Remark}
\newtheorem*{rem*}{Remark}

\begin{document}
\date{\today}
\title{On right conjugacy closed loops of twice prime order}
\author{Katharina Artic and Gerhard Hiss}


\address{K.A.: Lehrstuhl B f\"ur Mathematik, RWTH Aachen University,
52056 Aach\-en, Germany}
\address{G.H.: Lehrstuhl D f\"ur Mathematik, RWTH Aachen University,
52056 Aach\-en, Germany}

\email{katharina.artic@math.rwth-aachen.de}
\email{gerhard.hiss@math.rwth-aachen.de}


\begin{abstract}
The right conjugacy closed loops of order $2p$, where $p$ is an odd 
prime, are classified up to isomorphism.
\end{abstract}

\maketitle

\section{Introduction}\label{1}
\markright{INTRODUCTION}

A \textit{quasigroup}~$\mathcal{L}$ is a set with a binary operation 
$* : \mathcal{L} \times \mathcal{L} \rightarrow \mathcal{L}$, 
such that every equation $x*a = b$ or $a*x = b$ with $a, b \in \mathcal{L}$
has a unique solution~$x$. In this case, for every $a \in \mathcal{L}$, the 
\textit{right multiplication} $R_a: \mathcal{L} \rightarrow \mathcal{L}, 
x \mapsto x * a$ is a permutation of~$\mathcal{L}$ (and of course so is
every \textit{left multiplication}). A quasigroup is a \textit{loop}, if it 
contains an identity element. Thus a group is just a loop, in which the
operation is associative, and we will indeed view groups as loops.

In the following, we will only consider finite loops.
Let~$\mathcal{L}$ be a (finite) loop, whose identity element we denote by~$e$.
The \textit{right multiplication group of~$\mathcal{L}$} is the group
$G := \langle R_a \mid a \in \mathcal{L} \rangle$, a subgroup of the
symmetric group on~$\mathcal{L}$. Clearly,~$G$ acts faithfully and 
transitively on~$\mathcal{L}$ and~$R_e$ is the identity element of~$G$,
which we denote by~$1$. Let $H \leq G$ 
denote the stabiliser in~$G$
of $e \in \mathcal{L}$, and let $T := \{ R_a \mid a \in \mathcal{L} \}$. 
Then~$T$ is a 
transversal for $H^g \backslash G := \{ H^g x \mid x \in G \}$
for every $g \in G$, the identity
element of~$G$ is contained in~$T$, and $\langle T \rangle = G$. The 
triple~$(G,H,T)$ is called the \textit{envelope} of~$\mathcal{L}$,
a group theoretic object..

Conversely, starting from group theory, one defines
a \textit{loop folder} to be a triple $(G,H,T)$ of a finite group~$G$, a
subgroup $H \leq G$ and a subset $T \subseteq G$ with $1 \in T$, such 
that~$T$ is a transversal for $H^g \backslash G$ for every $g \in G$.
Given a loop folder $(G,H,T)$ one can construct a loop $(\mathcal{L},*)$
on the set $H\backslash G$ of right cosets 
of~$H$ in~$G$. However, the envelope of~$\mathcal{L}$ need not be equal
to~$(G,H,T)$. In contrast to the right multiplication group of~$\mathcal{L}$,
in general the group~$G$ will not act faithfully on~$\mathcal{L}$, and 
the transversal~$T$ will not generate~$G$. On the other hand,
it is not difficult to construct the envelope of~$\mathcal{L}$ 
from $(G,H,T)$. 

These results, as well as the notion of \textit{loop folder} and 
\textit{envelope of a loop} are contained in \cite[Section~$1$]{Asch}.
However, the connection between loops and their envelopes goes back to 
Baer~\cite{Baer}.

Let~$\mathcal{L}$ be a loop with envelop $(G,H,T)$.
We say that~$\mathcal{L}$ is \textit{right conjugacy closed}, or an
$\rcc$ loop, if $T = \{ R_a \mid a \in \mathcal{L} \}$ is
closed under conjugation by itself. Clearly, this is the case if and only if~$T$
is invariant under conjugation in $G = \langle T \rangle$; in other words,
if~$T$ is a union of conjugacy classes of~$G$. We shortly say 
that~$T$ is $G$-invariant in the following. Thus an $\rcc$ loop gives rise to 
a $G$-invariant transversal of~$H$, the stabilizer of~$e$ in~$G$. 
(A $G$-invariant transversal of a subgroup~$H$ of a group~$G$ is sometimes 
called a \textit{distinguished transversal} in the literature.)
On the group theoretic side, 
this leads to the notion of an $\rcc$ loop folder. This is a loop folder
$(G,H,T)$, where~$T$ is $G$-invariant. More definitions regarding loop
folders are given at the beginning of Section~\ref{3}.

It has been shown by Dr{\'a}pal~\cite{Drapal} that an $\rcc$ loop of prime 
order is a group. In this paper we determine all $\rcc$ loops of order $2p$,
where~$p$ is an odd prime. In order to achieve this, we first describe the
possible envelopes $(G,H,T)$ of such loops. Our approach is group theoretic.
In Sections~\ref{2} and~\ref{3} 
we show that if $(G,H,T)$ is an $\rcc$ loop folder such that~$G$ acts
faithfully on $H\backslash G$ and the index $|G\colon H|$ is the
product of two distinct primes, then~$G$ acts imprimitively on
$H\backslash G$ (Theorem~\ref{Imprimitive}). This result uses the classification 
of the finite simple groups and is based on the classification of finite
primitive permutation groups of squarefree degree by Li and Seress, and on
the determination of the minimal degrees of permutation representations of 
finite groups of Lie type by Patton, Cooperstein and Vasilyev. For the purpose
of our further investigation, it would suffice to enumerate the primitive
permutation groups of degree~$2p$ for odd primes~$p$; we are not aware of
any result in this direction which does not rely on the classification of the 
finite simple groups.

In Section~\ref{4}, we continue with some basic results on permutation groups 
of degree~$p$ and give a new proof of Dr{\'a}pal's theorem on $\rcc$
loops of prime order (Corollary~\ref{DrapalsTheorem}). 

Let $(G,H,T)$ be the envelope of an $\rcc$ loop of order $2p$, where~$p$ is an 
odd prime. Using Theorem~\ref{Imprimitive} mentioned above, we may now assume 
that there is a subgroup
$K \lneq G$ with $H \lneq K$, and also that one of the indices $|G\colon K|$
or $|K\colon H|$ is equal to~$2$, and the other index is equal to~$p$.
This configuration is analysed in Section~\ref{5} with elementary group 
theoretical methods. It turns out that there are three possible types for~$G$.
Firstly,~$G$
can be isomorphic to the wreath product $C_p \wr C_2$, where~$C_p$ and~$C_2$ 
denote (cyclic) groups of order~$p$ and~$2$, respectively. Secondly,~$G$ can be
isomorphic to a subgroup of $\Aff(1,p)$, the affine group over $\mathbb{F}_p$.
Thirdly,~$G$ can be isomorphic to a  group $K \times \langle a \rangle$, 
where~$K$ is an odd order subgroup of $\Aff(1,p)$ and~$a$ is an element of 
order~$2$ (Theorem~\ref{MainResult}). In particular,~$G$ is soluble.
Ultimately, our results rely on the classification of the finite simple groups.
One could avoid this by assuming from the outset that $G$ is soluble. This would
lead to exactly the same list of $\rcc$ loops of order $2p$, but of course 
without the guarantee to have found them all.

In Section~\ref{6}, we determine the number of isomorphism classes
of loops of order~$2p$ (Theorem~\ref{NumberOfAllLoops}).

Finally, Section~{7} introduces a series of examples of $\rcc$ loops of 
order $q^2 - 1$ and multiplication groups $\GL(2,q)$ 
(Proposition~\ref{GL2Examples}). For~$q = 4$, we obtain
a loop of order $3 \cdot 5$, whose multiplication group is not soluble.
These examples indicate that a generalisation of our results to $\rcc$ loops of
order $pq$ for distinct primes~$p$ and~$q$ could be substantially more difficult.

This is a good place to discuss some related results.
In~\cite[Theorem~A]{Stein}, Stein shows that if~$T$ is a conjugacy class in a
finite group and at the same time a tranversal to a subgroup, then
$\langle T \rangle$ is soluble. This result uses the classification of the
finite simple groups. Without the classification, but with the help of the
Odd Order Theorem, Cs\H{o}rg\H{o} and Niemenmaa in \cite{CsNie} obtain the
solubility of the full multiplication group of a loop under certain conditions
on the stabilizer of a point. Their paper contains further references for
results along this line. In \cite{CsDra}, Cs\H{o}rg\H{o} and Dr{\'a}pal
characterise left conjugacy closed loops inside the class of nilpotent loops
of nilpotency class two. In the same paper, these authors also determine the 
nilpotent left conjugacy closed loops of order $p^2$ for primes~$p$. In
\cite[Theorem~$4.15$]{Kunen}, Kunen shows that for each odd prime~$p$ there
is exactly one non-assciative conjugacy closed loop of order~$2p$, up to
isomorphism (a loop is conjugacy closed, if it is both left and right conjugacy
closed). Burn shows in \cite{Burn} that every Bol loop of order~$p^2$ or~$2p$
for a prime~$p$ is a group. Finally, in \cite[Theorem~$7.1$]{DalyVoj} Daly and 
Vojt\v{e}chovsk{\'y} determine
the number of nilpotent loops of order~$2p$, where again~$p$ is a prime,
up to isomorphism.

This paper builds upon the PhD thesis of the first author \cite{Artic}, written
under the direction of the second author and Alice Niemeyer. 
Theorem~\ref{Imprimitive} is contained in this thesis, but also a complete
classification of all $\rcc$ loops of order at most~$30$. These have been 
incorporated into the GAP package \textit{Loops} of Nagy and 
Vojt\v{e}chovsk{\'y}~\cite{NV15}. The classification of
the $\rcc$ loops of order~$2p$ is, to the best of our knowledge, new. The 
examples computed in~\cite{Artic}
were of considerable importance for confirming our theoretical results of 
Section~{6}.
The example of an $\rcc$ loop of order $15$ and multiplication group $\GL(2,4)$ 
contained in~\cite{Artic}, gave
rise to the series of examples constructed in Section~{7}.

Our group theoretical notation is standard. For example, we write $G'$ for
the commutator subgroup of the group~$G$. We do recall the notion of an
almost simple group and that of the core of a subgroup in the introductions 
to Section~\ref{2} and Section~\ref{3}, respectively.
As already indicated above, a cyclic group of order~$n$ is denoted by~$C_n$, 
and the symmetric and
alternating groups of degree~$n$ are denoted by~$S_n$ and~$A_n$, respectively.

\section{Primitive permutation groups of squarefree degree}\label{2}
\markright{PRIMITIVE GROUPS OF SQUAREFREE DEGREE}

We begin with a remark on the sizes of conjugacy classes in almost simple 
groups. Recall that a group~$G$ is \textit{almost simple}, if there is a 
non-abelian finite simple group~$S$ such that $S \leq G \leq \Aut(S)$ 
(where~$S$ is identified with the group of inner automorphisms of~$S$). In this 
context,~$S$ is called the \textit{socle} of~$G$.
\begin{rem}\label{SizesCC}
{\rm 
Let~$G$ be an almost simple group with socle~$S$. Denote by~$l$ the smallest 
index of any proper subgroup of~$S$. Since~$S$ is simple,~$l$ is a lower bound 
for the size of those non-trivial conjugacy classes of~$G$ lying in~$S$. Let 
$g \in G \setminus S$. Then we have
\begin{eqnarray*}
|G\colon C_G(g)| & = & \frac{|G|}{|SC_G(g)|} \cdot \frac{|S|}{|S \cap C_G(g)|} \\ 
& = &  \frac{|G|}{|SC_G(g)|} \cdot \frac{|S|}{|C_S(g)|}.
\end{eqnarray*}
Notice that if $C_{S}(g) = S$, the element~$g$ acts trivially on~$S$ which implies 
that $g = 1$. Hence we have $C_{S}(g) \lneq S$. Thus,~$l$ is a lower bound on 
the size of all non-trivial conjugacy classes of~$G$. 
}
\end{rem}

The following theorem combines some major results by Li and Seress on finite 
primitive permutation groups of squarefree degree, and by Patton, Cooperstein
and Vasilyev on the minimal degrees of permutation representations of finite 
groups of Lie type.

\begin{thm}
\label{SquareFree}
Let~$G$ be a finite primitive permutation group of degree~$n$ (i.e.\ $G$ acts 
faithfully and primitively on a set of~$n$ points). Suppose that~$n$ is 
square-free (i.e.\ $p^2 \nmid n$ for all primes~$p$). 
Then every non-trivial conjugacy class of~$G$ has at least~$n$ elements, or one
of the following holds:
\begin{itemize}
\item [{\rm (a)}] We have $n = p$ is a prime and~$G$ is isomorphic to a subgroup 
of $\Aff(1,p)$, 
\item [{\rm (b)}] We either have have $G=S_8$ and $n\in\{35,105\}$, or $G=J_1$ 
and $n = 2926$, or $G=\PGL(2,r)$ for an odd prime $r$ and $n=r(r+1)/2$,
\item [{\rm (c)}] or~$G$ is almost simple and $\soc(G)$ and~$n$ occur in 
{\rm Table~\ref{SocAlSimGr}}. There,~$r$ denotes a prime power. 
\end{itemize}
\end{thm}
\renewcommand{\baselinestretch}{1.3}
\begin{table}[t]
\begin{small}
\caption{\label{SocAlSimGr} Primitive groups of degree~$n$ which might have a
non-trivial conjugacy class of length less than~$n$}
\begin{center} 
\begin{tabular}{c | l | c | l }
& $\soc(G)$ & $n$ & Restrictions\\ \hline
(i) & $A_c$ & $\binom{c}{k}$ & $3\leq k\leq c-3$\\
(ii) & $A_{2a}$ & $\frac{1}{2}\binom{2a}{a}$ & $a\in\{6,9,10,12,36\}$\\
(iii) & $\PSL(m,r)$ & $\frac{\prod_{i=0}^{k-1}(r^{m-i}-1)} {\prod_{i=1}^k(r^i-1)}$ & $2\leq k\leq m-2$,\\
    & & & $(m,r)\notin\{(4,2),(5,2)\}$ \\
(iv) & $\PSL(m,r)$ & $\frac{\prod_{i=0}^{2k-1}(r^{m-i}-1)}{(\prod_{i=1}^k(r^i-1))^2}$ & $1\leq k<m/2$, $m\geq3$,\\
    & & & $(m,r)\neq(3,2)$ \\
(v) & $\PSL(2,r)$ & $\sqrt{r}(r+1)/2$ & $\sqrt{r}$ an odd prime,\\
    & & & $\soc(G)<G$, $r>9$ \\
(vi) & $\PSL(2,r)$ & $r(r^2-1)/24$ & $r$ a prime, $r \equiv \pm 3\ ({\rm mod}\ 8)$, \\
    & & & $r\notin\{5,11\}$\\
(vii) & $\PSL(2,r)$ & $r(r^2-1)/48$ & $r$ a prime, $r \equiv \pm 1\ ({\rm mod}\ 8)$, \\
    & & & $r\notin\{7,17,23\}$\\
(viii) & $\PSL(2,r)$ & $r(r^2-1)/120$ & $r$ a prime, $r \equiv \pm 1\ ({\rm mod}\ 10)$, \\
    & & & $r\notin\{11,19,29,31,41,59\}$\\
(ix) & $\PSU(4,r)$ & $(r^2 + 1)(r^3 + 1)$ & \\
(x) & $\PSp(2m,2)$ & $4^m - 1$ & $m \geq 3$ \\
(xi) & $\mathrm{PSp}(2m,r)$ & $\frac{(r^{2m}-1)(r^{2m-2}-1)}{(r^2-1)(r-1)}$ & $m\geq 3$ \\
(xii) & $\Omega(2m+1,r)$ & $\frac{(r^{2m}-1)(r^{2m-2}-1)}{(r^2-1)(r-1)}$ & $m\geq 3$ \\
(xiii) & $\mathrm{P}\Omega^-(2m,r)$ &$\frac{(r^m+1)(r^{2m-2}-1)(r^{m-2}-1)}{(r^2-1)(r-1)}$	& $m\geq 3$, $r$ even \\
(xiv) & $\mathrm{P}\Omega^-(2m,r)$ &$\frac{(r^m+1)(r^{2m-2}-1)(r^{2m-4}-1)(r^{m-3}-1)}{(r^3-1)(r^2-1)(r-1)}$& $m \equiv 0\ ({\rm mod}\ 4)$, $r$ even \\
(xv) & $\mathrm{P}\Omega^+(2m,2)$ & $(2^m - 1)(2^{m-1} + 1)$ & $m \geq 5$ odd \\
(xvi) & $\mathrm{P}\Omega^+(2m,r)$ & $\frac{(r^m-1)(r^{2m-2}-1)(r^{m-2}+1)} {(r^2-1)(r-1)}$ & $m \geq 3$, $r$ even \\
(xvii) & $\mathrm{P}\Omega^+(2m,r)$ & $\frac{(r^m-1)(r^{2m-2}-1)(r^{2m-4}-1)(r^{m-3}+1)}{(r^3-1)(r^2-1)(r-1)}$ & $m \equiv 3\ ({\rm mod}\ 4)$, $r$ even \\
(xviii) & $E_7(r)$ & $\frac{(r^{18}-1)(r^{14}-1)(r^4-r^2+1)}{(r^2-1)(r-1)}$ &  
\end{tabular}
\end{center}
\end{small}
\end{table}
\renewcommand{\baselinestretch}{1.1}

\begin{prf}
By \cite[Theorem~$1$]{LS} we either have that~$n$ is a prime and 
$G \leq \Aff(1,n)$ as in Case~(a), or~$G$ is almost simple and $S := \soc(G)$ as 
well as~$n$ appear in the paper~\cite{LS} by Li and Seress. 
	
The cases when~$S$ is isomorphic to an alternating group, are listed in 
\cite[Table~$1$]{LS}. If~$S$ is as in \cite[Table~$1$, Line~$1$]{LS}, then
$S = A_c$ and $n = \binom{c}{k}$ with $1 \leq k \leq c - 1$. For reasons of 
symmetry it suffices to consider the case $k \leq c/2$. 
Table~\ref{SmallAlternating} lists the size $s(G)$ of the smallest non-trivial 
conjugacy class of~$G$ for all~$G$ with $S \in \{A_5, A_6, A_7, A_8 \}$. This
table, easily compiled or verified with GAP~\cite{GAP4}, proves our claim for
$5 \leq c \leq 8$.
\begin{table}
\caption{\label{SmallAlternating} Smallest size of non-trivial conjugacy classes}
$$
\begin{array}{l|ccccccccccc}
G & A_5 & S_5 & A_6 & S_6 & A_6.2_2 & A_6.2_3 & \Aut( A_6 ) & A_7 & S_7 & A_8 & S_8 \\ \hline
s & 12 & 15 & 40 & 15 & 36 & 45 & 30 & 70 & 21 & 105 & 28 
\end{array}
$$
\end{table}
If $c \geq 9$, by \cite[Theorems~$5.2$A,B]{DM}, the subgroups 
of~$A_c$ or~$S_c$ which have an index less then $c(c-1)/2$ do not occur as 
centralizers of non-trivial elements. Hence the non-trivial conjugacy classes 
of~$A_c$ or~$S_c$ have at least $c(c-1)/2$ elements, proving our claim for 
$k \leq 2$. The case $3\leq k\leq c-3$ appears as Case (c)(i) in our statement.
	
In the remaining cases of \cite[Table~$1$]{LS}, a look at 
Table~\ref{SmallAlternating} shows that all non-trivial conjugacy classes of~$G$ 
have at least~$n$ elements except for
	\begin{itemize}
		\item $G=S_8$ and $n\in\{35,105\}$,
		\item $S = A_{2a}$ and $n = \binom{2a}{a}/2$ with $a \in 
                      \{6,9,10,12,36\}$.
	\end{itemize}
These cases appear as Case~(b) and Case~(c)(ii), respectively, in our statement.
	
The cases when~$S$ is a sporadic simple group are listed in \cite[Table~$2$]{LS}.
Using GAP, we only find the one exception listed in Case~(b).

In \cite[Table~$3$]{LS}, the case where~$S$ is a classical group are considered.
For some small parameter values, we have verified our claim directly with GAP. 
These cases are listed in the column headed \textit{Restrictions} of 
Table~\ref{SocAlSimGr}, and are not commented on any further below.
In the following, we refer to the line numbers of \cite[Table~$3$]{LS}.
Suppose that~$S$ is as in Line~$1$. Then $S = \PSL(m,r)$ and
	\begin{equation*}
		n=\prod_{i=0}^{k-1}(r^{m-i}-1) / \prod_{i=1}^k(r^i-1)
	\end{equation*}
with $1 \leq k < m$. For $k = 1$ or $k = m - 1$, we have $n = (r^m-1)/(r-1)$. 
If $(m,r) \in \{ (2,5), (2,7), (2,9), (2,11), (4, 2) \}$, a computation with
GAP shows that the non-trivial conjugacy classes of~$G$ have more than~$n$ 
elements. Otherwise,~$n$ is the smallest index of any proper subgroup of~$S$ 
by \cite[Table~$5.2$.A]{KL}. Applying Remark~\ref{SizesCC}, we see that the 
non-trivial conjugacy classes of~$G$ have at least~$n$ elements. The case 
$2 \leq k \leq m - 2$ is listed as Case~(c)(iii) in our statement. 
	
The case when~$S$ is as in Line~2, is listed as Case (c)(iv) in our statement. 
	
Suppose that~$S$ is as in Line~3 or~4. Then $S = \PSL(2,r)$ and 
$n = r(r \pm 1)/2$. Since~$n$ is squarefree, we have $r = 4$ and 
$n \in \{ 6, 10 \}$ or~$r$ is an odd prime. If $r = 4$, we have $S \cong A_5$,
and Table~\ref{SmallAlternating} proves our claim.
If~$r$ is an odd prime, then $\Aut(\PSL(2,r)) = \PGL(2,r)$ and hence 
$G = \PSL(2,r)$ or $G = \PGL(2,r)$. The conjugacy classes of these groups are 
well known. We find that only if $G = \PGL(2,r)$ and $n = r(r+1)/2$, there 
are non-trivial conjugacy classes of~$G$ with less than~$n$ elements. This case 
appears in Case~(b) in our statement
	
Suppose that~$S$ is as in Line~5. Then $S = \PSL(2,r)$ and 
$n = \sqrt{r}(r + 1)/2$. Since~$n$ is squarefree,~$r$ is the square of a 
prime number. The non-trivial conjugacy classes of $\PSL(2,r)$ have at 
least~$n$ elements. Hence $S \lneq G$. This case is listed as Case~(c)(v)
in our statement.
	
Suppose that~$S$ is as in one of the Lines 6, 7 or~8. Then $S = \PSL(2,r)$ and 
$n = r(r^2-1)/d$ with $d \in \{24,48,120\}$, and $r \equiv \pm 3\ ({\rm mod}\ 8)$ 
if $d=24$, respectively $r \equiv \pm 1\ ({\rm mod}\ 8)$ if $d=48$, respectively 
$r \equiv \pm 1\ ({\rm mod}\ 10)$ if $d=120$. In particular,~$r$ is odd. 
Since~$n$ is squarefree, $r = 9$ or~$r$ is an odd prime. If $r = 9$ we 
have $S \cong A_6$, and Table~\ref{SmallAlternating} proves our claim.
The cases where~$r$ is an odd prime, are listed as Case~(c)(vi) through 
Case~(c)(viii) in our statement.
	
Suppose that~$S$ is as in Line~9. Then $S = \PSU(m,r)$ with 
$$n = \frac{(r^m - (-1)^m)(r^{m-1} - (-1)^{m-1})}{r^2 - 1}.$$
If~$m = 2$, we have $S \cong \PSL(2,r)$ (see \cite[Theorem~$10.9$]{Taylor}) 
and $n = r + 1$, a case we have already 
considered above. For $m = 3$ and $r = 5$, our claim can be verified with GAP.
The case of $m = 4$ is listed as Case~(c)(ix) in our statement.
If $6 \mid m$ and $r = 2$, then $n = (2^m - 1)(2^{m-1} + 1)/3$ is not
squarefree, as $2^6 - 1$ divides $2^m - 1$ and~$3$ divides $2^{m-1} + 1$. 
In the remaining cases,~$n$ is the smallest index of any proper 
subgroup of~$S$ (see \cite[Table~$5.2$.A]{KL}). Thus by Remark \ref{SizesCC}, 
the non-trivial conjugacy classes of~$G$ have at least~$n$ elements.

Suppose that~$S$ is as in Line~10. Then $S = \PSp(2m,r)$ with $(m,r) \neq (2,2)$
and $n = (r^{2m} - 1)/(r - 1)$. 
(The case $(m,r) = (2,2)$ leads to $S = A_6$ and $n = 15$, which can be 
excluded by Table~\ref{SmallAlternating}.)
Again, we have already considered the case $m = 1$, where $S \cong \PSL(2,r)$
(see \cite[Theorem~$8.1$]{Taylor}). If $m = 2$ and $r = 3$,
then $n = 40$ is not squarefree.
The case $m \geq 3$ and $r = 2$ is listed
as Case~(c)(x) in our statement. In the remaining cases,~$n$ is the smallest 
index of any proper subgroup of~$S$ (see \cite[Table~$5.2$.A]{KL})  and
Remark \ref{SizesCC} proves our claim.

If~$S$ is as in Line~$12$ or~$13$, then $S \cong A_6$, and we are done with
Table~\ref{SmallAlternating}.

Suppose that~$S$ is as in Line~$14$. Then $S = \Omega(2m+1, r)$ and 
$n = (r^{2m} - 1)/(r - 1)$. We may assume that $m \geq 3$ and that~$r$ is odd, 
as otherwise $S \cong \PSp(2m,r)$ (see 
\cite[Theorems~$11.6$,~$11.9$, Corollary~$12.32$]{Taylor}), a case already 
considered. If $r = 3$,
then $n = (3^{2m} - 1)/2$ is not square free. In the other cases,~$n$ is the 
smallest index of any proper subgroup of~$S$ (see \cite[Table~$5.2$.A]{KL}),
and we are done as above.

Suppose that~$S$ is as in Line~$16$. Then,~$m$ is even and once more by 
\cite[Table~$5.2$.A]{KL} and Remark \ref{SizesCC} we obtain our claim. (This 
includes the case $m = 2$, where $S \cong \PSL(2,q^2)$ 
(see \cite[Corollary~$12.43$]{Taylor}) and $n = q^2 + 1$.)

Suppose that~$S$ is as in Line~$19$. Then $S = \mathrm{P}\Omega^+(2m,r)$ 
with $m \geq 3$ odd and
$$n = \frac{(r^m - 1)(r^{m-1} + 1)}{r-1}.$$
If $m = 3$, we have $S \cong \PSL(4,r)$ (see \cite[Corollary~$12.21$]{Taylor})
and $n = (q^2+1)(q^2+q+1)$. This case
is already contained in Case~(c)(iii) of our statement.
If $r \neq 2$ and $m \geq 5$, we conclude with \cite[Table~$5.2$.A]{KL} and 
Remark \ref{SizesCC}.
The case of $m \geq 5$ and $r = 2$ is included as Case (c)(xv) in our statement.

The remaining cases of \cite[Table~$3$]{LS} are listed as Cases (c)(xi) 
through (c)(xiv) and (c)(xvi) through (c)(xvii), respectively in our statement. 
	
Suppose that~$S$ is as in \cite[Table~$4$]{LS}, i.e.\ an exceptional group of 
Lie type.
In \cite{V1}, \cite{V2} and \cite{V3}, A.~V.~Vasilyev lists the smallest 
index~$l$ of any proper subgroup of the exceptional simple groups. By 
Remark~\ref{SizesCC}, the non-trivial conjugacy classes of~$G$ have at least~$l$ 
elements. We find $l = n$ except for
	\begin{itemize}
		\item	$S = G_2(4)$ and $n = 1365$. We verified our claim for 
the two almost simple groups with socle $G_2(4)$ with GAP.
		\item $S = E_7(r)$ with 
$n = (r^{18}-1)(r^{14}-1)(r^4-r^2+1)/((r^2-1)(r-1))$. This case is listed as 
Case~(c)(xviii) in our statement.
	\end{itemize}	
This completes our proof.
\end{prf}

\begin{rem}
\label{Table1Rem}
{\rm
For the purpose of this remark, let us call an \textit{example} a pair~$(G,n)$
of a primitive permutation group~$G$ of square-free degree~$n$ containing a 
non-trivial conjugacy class with less than~$n$ elements.

Now assume the hypotheses of Theorem~\ref{SquareFree}. Clearly, not all the 
instances $(G,n)$ listed there are examples.

If~$G$ is as in (a) of this theorem, then~$G$ has a conjugacy class of 
length~$p-1$. The symmetric group~$S_8$ has a conjugacy class of length~$28$,
and the sporadic simple group $J_1$ has a conjugacy class of length~$1463$.
Thus the groups in (a) and the first three instances of~(b) provide examples.
This fact also indicates that in order to enumerate all examples one will have
to use the classification of the finite simple groups.

The group $G = \PGL(2,r)$ has a conjugacy class of length $r(r-1)/2$ for every
odd prime~$r$. The group $\PGL(2,3)$ is isomorphic to the symmetric group~$S_4$,
which does not have a primitive permutation representation of degree~$6$.
Thus $(\PGL(2,r), r(r + 1)/2)$ is an example, if and only if $r \geq 5$ and 
$(r + 1)/2$ is square-free.

We expect that not many examples will arise from the pairs $(G,n)$ listed in 
Theorem~\ref{SquareFree}(c), but it would be a tedious task to enumerate all
of them. One approach could be to determine all subgroups of~$G$ of index
less than~$n$ and show that most of such subgroups are not centralizers of
elements. Still, one has to decide whether one of the remaining numbers~$n$ 
is indeed square-free. This will most certainly lead to difficult, if not 
intractible, number theoretical questions.
}
\end{rem}

In the lemma below we are going to make use of Zsigmondy primes, also known as 
primitive prime divisors. Let~$r$ and~$d$ be integers greater than~$1$. We call 
a prime~$\ell$ a Zsigmondy prime for $r^d - 1$, if $\ell$ divides $r^d - 1$, but 
not $r^i - 1$ for $1 \leq i < d$. A Zsigmondy prime for $r^d -1$ exists whenever 
$d > 2$ and $(r,d) \neq (2,6)$ (see \cite[Theorem IX.8.3]{HuBII}).

\begin{lem}\label{Grpspq}
Suppose that $n = pq$, where~$p$ and~$q$ are distinct primes, and that~$G$ is a 
finite primitive permutation group of degree~$n$ such that~$G$ has a non-trivial 
conjugacy class with less than~$n$ elements. Then one of the following holds: 
\begin{itemize}
\item [{\rm (a)}] We have $G\in\{A_7,S_7,S_8\}$ and $n=35$.
\item [{\rm (b)}] We have $G = \PGL(2,r)$ for an odd prime~$r$ and $n=r(r+1)/2$.
\item [{\rm (c)}] The group~$G$ is almost simple with 
$\PSL(2,r) = \soc(G) \lneq G$, where $\sqrt{r}$ is an odd prime, $r > 9$, and 
$n = \sqrt{r}(r+1)/2$.
\item [{\rm (d)}] We have $G \in \{\PSL(2,13),\PGL(2,13)\}$ and $n=91$.
\item [{\rm (e)}] We have $G \in \{\PSL(2,61),\PGL(2,61)\}$ and $n=1891$.
\item [{\rm (f)}] We have $S = \mathrm{P}\Omega^+(2m,2)$, for $m \geq 3$
and $ n = (2^m - 1)(2^{m-1} + 1)$.
\end{itemize}
\end{lem}
\begin{prf}
We have to exclude those integers~$n$ in Theorem~\ref{SquareFree} which are not 
the product of two different primes. From Cases~(a) and~(b) of 
Theorem \ref{SquareFree}, we obtain (part of) Case~(a) and Case~(b) of our lemma.
	
So suppose that~$G$ is almost simple and that $S := \soc(G)$ occurs in 
Table~\ref{SocAlSimGr}.
In Case~(i) we have $n = \binom{c}{k}$ with $k \geq 3$. For reasons of symmetry 
it suffices to consider $3 \leq k \leq c/2$. By \cite[Theorem~$7$]{Ple}, the 
total number (counting multiplicities) of prime factors of the binomial 
coefficient $\binom{c}{k}$ is greater than or equal to the total number of prime 
factors of~$c$, with equality only if 
$(c , k)=(8,4)$. 
Thus $\binom{c}{k}$ is the product of two different primes only if~$c$ is a 
prime. Consider the case $k = 3$ first. Then
\begin{equation*}
n = c \cdot \frac{(c-1)(c-2)}{6}.
\end{equation*}
If $c = 7$ we have $n = 35$. This is recorded in Case (a) of our lemma. 
So assume that $c \geq 11$. Then $n = \binom{c}{3}$ has at least three different 
prime factors. But then by \cite[Theorem~$3$]{Ple}, the binomial coefficient 
$\binom{c}{k}$ with $k > 3$ also has at least three different prime factors.

In Cases~(ii),~(ix) through~(xiv) and~(xvi) through~(xviii) of 
Table~\ref{SocAlSimGr}, the degree~$n$ is clearly not the product of two 
different primes.
	
In Case (iii) of Table~\ref{SocAlSimGr} we have $\soc(G)=\PSL(m,r)$ and
\begin{equation*}
n = \prod_{i=0}^{k-1}(r^{m-i}-1)/\prod_{i=1}^k(r^i-1),
\end{equation*}
with $2 \leq k < m$. 
For reasons of symmetry it suffices to consider the integers $k$ with 
$1 \leq k \leq m/2$. Suppose first that $m \geq 5$ and $k \geq 3$. Consider the 
terms $(r^m - 1)$, $(r^{m-1} - 1)$ and $(r^{m-2} - 1)$. They occur only in the 
numerator and not in the denominator of~$n$ and, by Zsigmondy's theorem, have 
pairwise distinct primitive prime divisors~$r_1$,~$r_2$ and~$r_3$ (which 
divide~$n$) unless one of the pairs $(m,r)$, $(m-1,r)$ or $(m-2,r)$ is equal
to~$(6,2)$. But in these cases, i.e.\ if $m \in \{6,7,8\}$ and $r = 2$, we just 
compute that~$n$ is not the product of two different primes for all 
$3 \leq k \leq m/2$. If $m \geq 4$ and $k = 2$, then
\begin{equation*}
n = \frac{(r^m-1)(r^{m-1}-1)}{(r^2-1)(r-1)},
\end{equation*}
which is the product of two different primes if and only if $m \in \{4,5\}$ and 
$r = 2$. But these cases have already been excluded.
	
In Case~(iv) of Table~\ref{SocAlSimGr} we have $S = \PSL(m,r)$ and
\begin{equation*}
n = \prod_{i=0}^{2k-1}(r^{m-i}-1)/(\prod_{i=1}^k(r^i-1))^2,
\end{equation*}
with $m \geq 3$ and $1 \leq k < m/2$. Suppose first we have $m \geq 5$ and 
$k \geq 2$. Then, again, the terms $(r^m - 1)$, $(r^{m-1} - 1)$ and 
$(r^{m-2} - 1)$ occur only in the numerator and not in the denominator of~$n$ 
and have pairwise distinct primitive prime divisors~$r_1$,~$r_2$ and~$r_3$ 
(which divides~$n$) unless one of the pairs $(m,r)$, $(m-1,r)$ or $(m-2,r)$ is 
equal to~$(6,2)$. But in these cases, i.e.\ if $m \in \{6,7,8\}$ and $r = 2$, we 
just compute that~$n$ is not the product of two different primes for all 
$2 \leq k<m/2$. If~$m$ is arbitrary and $k = 1$ then
\begin{equation*}
n = \frac{(r^m-1)(r^{m-1}-1)}{(r-1)(r-1)}
\end{equation*}
which is the product of two different primes if and only if $(m,r) = (3,2)$. But 
this case has already been excluded.
	
Case~(v) of Table~\ref{SocAlSimGr} is listed as Case~(c) in our lemma.
	
In Cases~(vi) through~(viii) of Table~\ref{SocAlSimGr}, we have 
$n = r(r^2-1)/d$ with $d\in\{24,48,120\}$. Clearly,~$n$ is not the product 
of two different primes if 
$r > d+1$. For $r \leq d+1$ and~$r$ not equal to one of the primes excluded in
Table~\ref{SocAlSimGr}, we have that~$n$ is the product of two different primes 
if and only if $(r,d) = (13,24)$ or $(r,d) = (61,120)$. We have $S = \PSL(2,r)$, 
and as~$r$ is a prime, $G = \PSL(2,r)$ or $G = \PGL(2,r)$. These cases 
are listed as Case~(d) and Case~(e) of our lemma.

Finally, Case~(xv) of Table~\ref{SocAlSimGr} is Case~(f) of the lemma.
\end{prf}

\section{The action of the right multiplication groups of rcc loops of \label{3}
twice prime order}
\markright{ACTION OF RIGHT MULTIPLICATION GROUPS OF RCC LOOPS}

The results of the previous section are now applied to $\rcc$ loops whose order 
is the product of two distinct primes. Recall the notion of the envelope of a 
loop as introduced in the second paragraph of Section~\ref{1} (which follows 
\cite[p.~$100$]{Asch}). Recall also that a loop folder is a triple $(G,H,T)$,
where~$G$ is a finite group, $H$ is a subgroup of~$G$ and~$T$ is a transversal,
with $1 \in T$, for all coset spaces $H^g\backslash G$ with $g \in G$ (see
the third paragraph of Section~\ref{1} and \cite[p.~$101$]{Asch}). A loop folder 
$(G,H,T)$ is an $\rcc$ loop folder, if~$T$ is invariant under conjugation by~$G$  
(see also Section~\ref{1}). We will now introduce further notation, 
although only needed in later sections.

Let $(G,H,T)$ be a loop folder.
By definition, the \textit{order} of $(G,H,T)$ is the size of~$T$.
We say that $(G,H,T)$ is \textit{faithful}, if~$G$ acts faithfully on
$H\backslash G$. This is the case if and only if the core of~$H$ in~$G$ is 
trivial. Recall that the smallest normal subgroup~$C$ of~$G$ containd in~$H$ is 
called the \textit{core of~$H$ in~$G$}. Thus~$C$ is the intersection of all the 
$G$-conjugates of~$H$ in~$G$, i.e.\ $C := \cap_{g \in G} H^g$. The core of~$H$ 
in~$G$ is equal to the kernel of the permutation representation of~$G$ on the 
(right or left) cosets of~$H$. Clearly, the envelope of a loop is a faithful 
loop folder.

Here is the main result of this section. It is only used in the setup of
Subsection~\ref{Generalities}, and nowhere else in this paper.

\begin{thm}
\label{Imprimitive}
Let $(G,H,T)$ denote the envelope of an $\rcc$ loop of order $n = pq$, where~$p$ 
and~$q$ are distinct primes. Then~$G$ acts imprimitively on~$H\backslash G$.
\end{thm}
\begin{prf}
Let $n = pq = |H\backslash G|$. Suppose that~$G$ acts primitively 
on~$H\backslash G$. Since $|T| = n$ and~$T$ is a union of conjugacy classes one
of which is the trivial class,~$G$ has a non-trivial conjugacy class with less 
than~$n$ elements. Hence~$G$ is one of the groups of Lemma~\ref{Grpspq}.

In Cases~(a),~(b),~(d) and~(e) of Lemma~\ref{Grpspq}, the concerned groups have 
at most two non-trivial conjugacy classes with less than~$n$ elements. 
Elementary combinatorics shows that in these cases there is no union of 
conjugacy classes~$T$ with $|T| = n$.

Suppose that~$G$ is as in Case~(c) of Lemma \ref{Grpspq}. Then~$G$ is almost 
simple with $S := \soc(G) = \PSL(2,r)$, where $\sqrt{r}$ an odd prime, $r > 9$, 
and $n = \sqrt{r}(r+1)/2$. Moreover, $S \lneq G$. The subgroups of 
$S = \PSL(2,r)$, are classified in Dickson's Theorem; see 
\cite[Hauptsatz~II.8.27]{HuI}. This shows that if $r \geq 17$, then only the 
maximal subgroups of index $r + 1$ have an index less then 
$n = \sqrt{r}(r+1)/2$. Consider the non-trivial conjugacy classes of~$G$. Those 
which contains elements of~$S$ have at least~$n$ elements as we already 
mentioned in the proof of Theorem~\ref{SquareFree}. 

Let $g\in G \setminus S$. By Remark \ref{SizesCC} we have
\begin{equation*}
|G\colon C_G(g)| = a \cdot |S\colon C_{S}(g)|
\end{equation*}
for some positive integer~$a$. Hence $|G\colon C_G(g)| \geq n$ except if 
$C_{S}(g) \leq M$, where~$M$ is a maximal subgroup of~$S$ with 
$|S\colon M| = r + 1$. In this case $r + 1 \mid |S\colon C_{S}(g)|$. Hence, if 
$|G\colon C_G(g)|$ is less than~$n$, it is a multiple of $r+1$. Thus, a union~$T$ 
of conjugacy classes of sizes less than~$n$ with $1 \in T$ has a size congruent 
to~$1$ modulo~$r+1$ and is therefore not divisible by the prime $(r+1)/2$. 
Therefore, it is not possible to have $|T| = n$.

Finally, assume that~$G$ is as in Case~(f) of Lemma~\ref{Grpspq}. Then~$G$ is 
almost simple with $S := \soc(G) = \mathrm{P}\Omega^+(2m,2)$, where $m \geq 3$ 
and $n = (2^m - 1)(2^{m-1} + 1)$. In order for~$n$ to be the product of 
two distinct primes, it is necessary that $p := 2^{m-1} + 1$ is a Fermat prime
and $q := 2^m - 1$ is a Mersenne prime. In particular,~$m$ is a Fermat prime. 
By \cite[Table~$5.2$.A]{KL}, the smallest index of any proper subgroup of~$S$ 
equals $2^{m-1}(2^m - 1)$. Remark~\ref{SizesCC} shows that any nontrivial 
conjugacy class of~$G$ has at least $2^{m-1}(2^m - 1)$ elements. As twice this 
number is greater than~$n$, we conclude that~$T$ is the union of two conjugacy 
classes of~$G$, one of which has length $2(2^{2m-2} + 2^{m-2} - 1)$. 
We have $m \geq 5$ and $S = \SO^+(2m,2)$ (in the notation of 
\cite[p.~$160$]{Taylor}). In particular, $G \in \{ \SO^+(2m,2), O^+(2m,2) \}$. 
For the order of~$G$ see \cite[p.~$141$,~$165$]{Taylor}.

Let~$\ell$ be a Zsigmondy prime for $2^{2m-2} - 1$. As~$\ell$ does not divide 
$2^{m-1} - 1$, and as $p = 2^{m-1} + 1$, we conlcude that $\ell = p$. Let~$t$ be 
a nontrivial element in~$T$. Now~$p$ does not divide $2(2^{2m-2} + 2^{m-2} - 1)
= |G\colon C_G(t)|$, and so there is $g \in C_G(t)$ with $|g| = p$. Let~$V$ 
denote the natural $2m$-dimensional $\mathbb{F}_2$-vector space of~$G$, equipped
with the quadratic form~$Q$ defining~$G$. Since~$p$ is a Zsigmondy prime for 
$2^{2m-2} - 1$, it follows that~$g$ acts irreducibly on some 
$(2m - 2)$-dimensional subspace~$V_0$ of~$V$. As the dimension of~$V_0$ is 
larger than~$1$, either~$V_0$ is totally singular or non-degenerate with respect 
to~$Q$ (for these notions see \cite[p.~$56$]{Taylor}). The maximal dimension of a 
totally singular subspace of~$V$ equals~$m$, and thus~$V_0$ is in fact 
non-degenerate. It follows that~$g$ fixes $V_1 := V_0^\perp$. In particular, 
$g \in O(V_0) \times O(V_1)$, where $O(V_i)$ denotes the orthogonal group with 
respect to the restriction~$Q_i$ of~$Q$ to~$V_i$, $i = 0, 1$. We may thus write
$g = g_0 \oplus g_1$, whith~$g_i$ the restriction of~$g$ to~$V_i$, 
$i = 0, 1$. Now~$O(V_0)$ contains a cyclic, irreducible subgroup, and thus the 
Witt index of~$Q_0$ equals~$m - 2$ by \cite[Satz~3c)]{Hu}. Hence~$V_1$ does not 
contain any non-trivial singular vector with respect to~$Q_1$, and thus 
$O(V_1) \cong S_3$, the symmetric group on three letters (see 
\cite[Theorem~$11.4$]{Taylor}). As $|g| = p \geq 17$, we conclude that~$g$ acts 
trivially on~$V_1$, i.e.\ $g_1 = 1$ and~$V_1$ is the fixed space of~$g$. It 
follows that $C_G( g )$ fixes~$V_1$ and $V_0 = V_1^\perp$, and thus 
$C_G( g ) \leq C_{O(V_0)} (g_0 ) \times O(V_1)$. As~$g_0$ acts irreducibly 
on~$V_0$, its centralizer in~$O(V_0)$ is cyclic and irreducible, and thus equals 
$\langle g_0 \rangle$, again by \cite[Satz~3c)]{Hu}. In particular, 
$t \in C_G(g) \leq \langle g_0 \rangle \times O(V_1)$. If $p \mid |t|$ or
$3 \mid |t|$, then $C_G( t ) \leq O(V_0) \times O(V_1)$. Otherwise, $|t| = 2$
and~$t$ has a $(2m - 1)$-dimensional fixed space~$V'$ on~$V$ and
$C_G(t) \leq O(V')$. In any case, the $2$-part of $|C_G(t)|$ is less than
$2^{(m-1)^2}$, whereas the $2$-part of~$|G|$ equals $2^{m(m-1) + 1}$ (see
\cite[p.~$141$]{Taylor}). In particular, the $2$-part of $|G\colon C_G(t)|$
is larger than~$2$, a contradiction.
\end{prf}

We end this section with two general results on loop folders with certain
invariance properties. The first will be used in an extension of a theorem of 
Dr{\'a}pal~\cite{Drapal}. 

\begin{lem}
\label{Normalizers}
Let $(G,H,T)$ be a loop folder such that~$T$ is invariant under conjugation
by~$H$. Suppose the $t \in T$ is such that
$Ht$ is also a left $H$-coset in~$G$, i.e.\ there exist $g \in G$ with
$Ht = gH$ (this is the case in particular if~$t$ normalizes~$H$). Then
$[t,H] = 1$.
\end{lem}
\begin{prf}
Let $h \in H$. Then
$$Ht = gH = gHh = Hth = Hh^{-1}th.$$
This implies $t = h^{-1}th$, as $t, h^{-1}th \in T$.
\end{prf}

\begin{lem}
\label{LengthsOfConjugacyClasses}
Let $(G,H,T)$ denote an $\rcc$ loop folder. Let $K \leq G$ such that
$HG' \leq K$. Then
$$|G\colon C_G(t)| \leq |K\colon H|\quad\text{for all\ } t \in T$$
(i.e.\ the length of the conjugacy classes of the elements in~$T$ are
bounded above by $|K\colon H|$).
\end{lem}
\begin{prf}
Let $g \in G$. Then the right coset $Kg$ is a union of exactly $|K\colon H|$
right cosets of~$H$ in~$G$. Thus $|Kg \cap T| = |K\colon H|$.

Now let $t \in T$ and $x \in G$. Then $t^xt^{-1} \in G' \leq K$. It follows
that $Kt^x = Kt$ and thus $t^x \in Kt$ for all $x \in G$. As $t^x \in T$ 
for all $x \in G$ by assumption, we conclude that $|G\colon C_G(t)| = 
|\{ t^x \mid x \in G \}| \leq |K\colon H|$.
\end{prf}

\section{Right conjugacy closed loops of prime order}\label{4}
\markright{RIGHT CONJUGACY CLOSED LOOPS OF PRIME ORDER}

In this section we give a new proof of a theorem of Dr{\'a}pal~\cite{Drapal}
which states that left conjugacy closed loops of prime order are groups.
We prove the analogue for $\rcc$ loops, but as the opposite 
loop of a left conjugacy closed loop is an \rcc loop, our version is equivalent 
to Dr{\'a}pal's result. Recall the notions related to loop folders summarized
at the beginning of Section~{3}.

We begin with an easy lemma.

\begin{lem}
\label{CentralizersInSp}
Let~$p$ be a prime and let $G \leq S_p$ with $p \mid |G|$.
Then the following statements hold for every $1 \neq g \in G$.

{\rm (a)} If $p \nmid |g|$, then $p \mid |G\colon C_G(g)|$.

{\rm (b)} If $p \mid |g|$, then $|C_G(g)| = p$.

{\rm (c)} Suppose that $|G\colon C_G(g)| < p$.
Then~$G$ has a 
unique Sylow $p$-subgroup~$P$ and $g \in P$. Moreover, if $G \neq P$, then~$G$
is a Frobenius group with kernel~$P$ and a Frobenius complement of order $r$ 
dividing~$p - 1$. In this case,~$P$ is the Frattini subgroup of~$G$.

{\rm (d)} We have $|G\colon C_G(g)| \neq 2(p - 1)$.
\end{lem}
\begin{prf}
In view of the cycle decomposition of~$g$, the first two parts are trivial. 
So let us assume that $|G\colon C_G(g)| < p$ or that $|G\colon C_G(g)| = 
2(p - 1)$. By~(a) and~(b) we have $|C_G(g)| = p$, and in particular $P := 
\langle g \rangle$ is a Sylow $p$-subgroup of~$G$. Under the hypothesis of~(c) 
we get $|G\colon N_G(P)| < p$, and under the hypothesis of~(d) we get
$|G\colon N_G(P)| \mid 2(p - 1)$. In each case Sylow's theorems imply 
$P \unlhd G$. Now~(d) and the last two statements of~(c) follow from the 
fact that~$G$ embeds into $N_{S_p}(P)$, which is isomorphic to $\Aff(1,p)$.
\end{prf}

\begin{cor}[Dr{\'a}pal~\cite{Drapal}]
\label{DrapalsTheorem}
Let~$p$ be a prime and let $(G,H,T)$ denote the envelope of an $\rcc$ 
loop~$\mathcal{L}$ of order~$p$. Then $H = 1$, i.e.\ $\mathcal{L}$ is a group
(isomorphic to~$G$).
\end{cor}
\begin{prf}
We may assume that $G \leq S_p$ and we have $p \mid |G|$. Now $T = \{ 1 \} \cup 
T'$ with $T' := T \setminus \{ 1 \}$. By assumption, $T'$ is a union of 
conjugacy classes of~$G$ of lengths at most $p - 1$. It follows from
Lemma~\ref{CentralizersInSp}(c) that~$G$ has a unique Sylow $p$-subgroup~$P$
and that $T \subseteq P$.
Hence $G = \langle T \rangle = P$, i.e.\ $H = 1$.
\end{prf}

We will also need the following generalization of 
Corollary~\ref{DrapalsTheorem}.

\begin{prp}
\label{PrimeOrderLoopFolders}
Let~$p$ be a prime and let $(G,H,T)$ be an $\rcc$ loop folder of order~$p$
with $\langle T \rangle = G$. Then~$G$ is abelian.
\end{prp}
\begin{prf}
Let 
$$N := \bigcap_{g \in G} H^g$$
denote the kernel of the action of~$G$ on $H\backslash G$ and put
$\bar{G} := G/N$, $\bar{H} := H/N$ and $\bar{T} := \{ Nt \mid t \in T \} 
\subseteq \bar{G}$. Then $(\bar{G}, \bar{H}, \bar{T})$ is a faithful $\rcc$ 
loop folder of order~$p$ with $\langle \bar{T} \rangle = \bar{G}$.
Thus $(\bar{G}, \bar{H}, \bar{T})$ is the envelope of an $\rcc$ loop of
order~$p$ (see \cite[$1.7.(4)$]{Asch}). By the result of Dr{\'a}pal (see 
Corollary~\ref{DrapalsTheorem}), 
such a loop is a group. It follows that $\bar{T} = \bar{G}$ and $\bar{H} = 1$. 
In particular, $H = N$ is a normal subgroup of~$G$ of index~$p$.

To show that~$G$ is abelian, let~$t \in T$.  By Lemma~\ref{Normalizers}
we have $[t,H] = 1$, as $H \unlhd G$. It follows that $\langle H, t \rangle
\leq C_G(t)$. If $t \neq 1$, we have $\langle H, t \rangle = G$, as~$H$ is
of index~$p$ in~$G$. Hence $t \in Z(G)$ for all $t \in T$. The claim follows
from $\langle T \rangle = G$.
\end{prf}

Notice that the above proposition is a generalization of Dr{\'a}pal's theorem
(see Corollary~\ref{DrapalsTheorem}); indeed, if $(G,H,T)$ is the envelope of
a loop, then $G = \langle T \rangle$, and if, moreover,~$G$ is abelian, then
$H = 1$, as the core of~$H$ in~$G$ is trivial.

\section{The right multiplication groups of rcc loops 
of twice prime order}\label{5}
\markright{RIGHT MULTIPLICATION GROUPS OF RCC LOOPS}

We refer the reader to the introduction of Section~\ref{3} for the notions
related to loop folders.

\subsection{Generalities}\label{Generalities}
Let~$p$ and~$q$ be distinct primes and let $(G,H,T)$ denote the envelope of an 
$\rcc$ loop of order $pq$. This implies in particular that $G$ acts faithfully 
on $H\backslash G$, i.e.\ the core of~$H$ in~$G$ is trivial. It is at this 
stage, and only here, where we impose an important consequence of 
Theorem~\ref{Imprimitive}. This states that~$G$ acts imprimitively on 
$H\backslash G$, and hence~$H$ is not a maximal subgroup of~$G$.
We let $K \lneq G$ such that $H \lneq K$. We choose 
notation such that $|G\colon K| = q$ and $|K\colon H| = p$.

Put $T_1 := T \cap K$, $K_1 := \langle T_1 \rangle \leq K$ and $H_1 := 
H \cap K_1 \leq H$.

We collect first properties.

\begin{lem}
\label{FirstProperties}
Let the notation be as above. Then $(K_1,H_1,T_1)$ is an $\rcc$ loop folder of
order~$p$ with~$K_1$ abelian. Also, $K_1 \unlhd K$ and $K = HK_1$. Finally,
$H_1 \unlhd K$.
\end{lem}
\begin{prf}
Clearly,~$K$ is the disjoint union of the cosets~$Ht$ for $t \in T_1$. Thus
$|T_1| = p$ and $K_1$ is the disjoint union of the cosets $K_1t$ for 
$t \in T_1$. As~$T_1$ is invariant under conjugation in~$K$, the first statement 
follows. The second statement follows from 
Proposition~\ref{PrimeOrderLoopFolders}, and the next two are obvious. The
last statement follows from $K = HK_1$ and the fact that $H_1 \unlhd H$ and 
that~$K_1$ is abelian.
\end{prf}

\subsection{The case $q = 2$.} \label{CaseQEqal2}Let us assume throughout this 
subsection that 
$q = 2$. Then $K \unlhd G$. Moreover, $K_1 \unlhd G$, as $K_1^g = 
\langle (T \cap K)^g \rangle = \langle T \cap K \rangle = K_1$ for all 
$g \in G$.

\begin{lem}
\label{NormalSubgroupsOfK}
Let $L \leq H$ with $L \unlhd K$. Then $L \cap L^a = 1$ and $LL^a =
L \times L^a \unlhd G$ for all $a \in G \setminus K$.
\end{lem}
\begin{prf}
Let~$a \in G \setminus K$. Clearly, $L \cap L^a$ and $LL^a$ are normal subgroups 
of~$G$, as $a^2 \in K$, and $L \unlhd K$. Thus $L \cap L^a = 1$ since 
$L \cap L^a \leq H$ and the core of~$H$ in~$G$ is trivial. As $L^a \unlhd K$, 
the product $LL^a$ is direct.
\end{prf}

We record two consequences which will be used later on.

\begin{cor}
\label{StructureOfK1}
We have $H_1 \in \{ 1, p \}$ and $K_1$ is elementary abelian of order~$p$ 
or~$p^2$.
\end{cor}
\begin{prf}
If~$H_1$ is trivial,~$K_1$ has order~$p$ by Lemma~\ref{FirstProperties}. Suppose 
that~$H_1$ is nontrivial and let $a \in T \setminus K$. Then 
$H_1 \cap H_1^a = 1$ by Lemma~\ref{NormalSubgroupsOfK}. As $K_1 \unlhd G$, we 
have $H_1H_1^a \leq K_1$.  It follows that $|H_1|^2$ divides $|K_1| = p |H_1|$. 
This implies $|H_1| = p$ and $K_1 = H_1H_1^a$, yielding our claim.
\end{prf}

\begin{cor}
\label{HNormalInK}
Suppose that $1 \neq H \unlhd K$. Then $|H| = p$ and there is an involution 
$a \in G \setminus K$ such that $K = H \times H^a$. In particular,~$G$ is
isomorphic to the wreath product $C_p \wr C_2$.
\end{cor}
\begin{prf}
By Lemma~\ref{NormalSubgroupsOfK} we have $H \cap H^a = 1$ and $HH^a \leq K$
for every $a \in G \setminus K$. As $|K| = p|H|$, this implies that $|H| = p$
and $K = H \times H^a$. It also follows that the involutions in~$G$ are 
contained in $G \setminus K$ and thus $G \cong C_p \wr C_2$.
\end{prf}

We now distinguish two cases.

\subsubsection{Case~$1$} \label{Case1Section}
Assume that $H \neq 1$ and that $[s, K_1] = 1$ for all $s \in T \setminus K$.
\begin{prp}
\label{Case1}
Under the assumptions of~{\rm \ref{Case1Section}}, we have $H_1 = 1$
and $K_1 \leq Z(G)$. Moreover,~$G$ is isomorphic to the wreath product 
$C_p \wr C_2$.
\end{prp}
\begin{prf}
The fact that~$K_1$ is abelian and our hypothesis imply that 
$T \subseteq C_G( K_1 )$, and hence $G = \langle T \rangle
\leq C_G(K_1)$, i.e.\ $K_1 \leq Z(G)$. Thus $H_1 \unlhd G$, which implies
$H_1 = 1$, as $H_1 \leq H$ and the core of~$H$ in~$G$ is trivial. Now
$K = HK_1$ by Lemma~\ref{FirstProperties}, and thus $H \unlhd K$. The claim 
follows from Corollary~\ref{HNormalInK}.
\end{prf}

\subsubsection{Case~$2$} \label{Case2Section}
Assume that there is $s \in T \setminus K$ 
with $[s,K_1] \neq 1$. In this case we put $Z := K_1 \cap C_G(s)$. Also, we
let $C := \cap_{k \in K} H^k \unlhd K$ denote the kernel of the action of~$K$
on the cosets of~$H$ in~$K$.

\begin{lem}
\label{Case2Properties}
Under the assumptions and with the notation of~{\rm \ref{Case2Section}}, 
the following statements hold.

{\rm (a)} We have $|G\colon C_G(s)| = p$ and $|Z| = |H_1|$.

{\rm (b)} We have $G = C_G(s)K_1$ and $Z \leq Z(G)$.

{\rm (c)} The centralizer $C_G(s)$ is abelian and $C_G(s)/Z$ is cyclic.
\end{lem}
\begin{prf}
By assumption, $K_1 \not\leq C_G(s)$. Corollary~\ref{StructureOfK1} implies
that $|Z| \in \{1, p \}$. By Lemma~\ref{FirstProperties} we have $|T_1| = p$,
which implies that $|T \setminus K| = p$ (recall that $q = 2$ and thus $|T| 
= 2p$). As $T \setminus K$ is a union of
conjugacy classes of~$G$, we have $|G\colon C_G(s)| \leq p$, i.e.\ 
$|C_G(s)| \geq |G|/p$. 

Now if $H_1 = 1$, i.e.\ $|K_1| = p$, we also have
$|Z| = 1$ and $G = C_G(s)K_1$. Thus all statements of (a) and (b) hold in 
this case. 

Now assume that $|H_1| = p$, i.e.\ $K_1$ is elementary abelian of 
order $p^2$. Then $|Z| \leq p$ and
$$|G| \geq |C_G(s)K_1| \geq \frac{|C_G(s)||K_1|}{|Z|} \geq \frac{|G|}{p} \cdot 
\frac{p^2}{|Z|} \geq |G|,$$
and we must have equality everywhere in the above chain of inequalities.
This implies $|Z| = p$ and 
$|C_G(s)K_1| =|G|$, again yielding all the claims of~(a) and the first
claim of~(b).

In any case, the set $T \setminus K$ is a conjugacy class of~$G$, consisting
of the elements $\{s^k \mid k \in K_1 \}$. Write $\bar{\ } : G \rightarrow 
\bar{G} := G/K_1$ for the canonical epimorphism. We have $\bar{G} = 
\langle \bar{T} \rangle = \langle \bar{s} \rangle$ as $T_1 \subseteq K_1$.
The natural isomorphism $\bar{G} \rightarrow C_G(s)/Z$ maps $\bar{s}$
to $Zs \in C_G(s)/Z$. Thus $C_G(s)/Z = \langle Zs \rangle$ and $C_G(s) =
\langle Z, s \rangle$. In particular, $C_G(s)/Z$ is cyclic and $C_G(s)$ is 
abelian, as $Z \leq C_G(s)$. Now $G = C_G(s)K_1$ and $Z = K_1 \cap C_G(s)$
imply that $Z \leq Z(G)$.
\end{prf}

The previous result implies in particular that~$G$ is soluble. Indeed,~$K_1$
is a normal subgroup of~$G$, as we have remarked at the beginning 
of~\ref{Case2Section}. Now~$K_1$ is abelian by Corollary~\ref{StructureOfK1},
and $G/K_1 = C_G(s)K_1/K_1 \cong C_G(s)/Z$ by Part~(b) of the lemma above and
by the definition of~$Z$. By Part~(c) of the lemma, $C_G(s)/Z$ is cyclic,
hence~$G$ is soluble.

\begin{cor}
\label{Case2Special}
Let the assumptions and notation be as in~{\rm \ref{Case2Section}}.
If $H_1 \neq 1$, then $H \unlhd K$.
\end{cor}
\begin{prf}
Suppose that $H_1 \neq 1$. Then $|Z| = p$ by Corollary~\ref{StructureOfK1}
and Lemma~\ref{Case2Properties}(a).
Now $Z \leq Z(G)$ by Lemma~\ref{Case2Properties}(b), and thus $Z \cap H = 1$,
since the core of~$H$ in~$G$ is trivial. It follows that $K = H \times Z$, as 
$|K\colon H| = p$ and $Z \leq K_1 \leq K$.

Now~$K/C$ is a soluble permutation group on~$p$
points. It follows from a theorem of Galois (see \cite[Satz~II.3.6]{HuI})
that $K/C$ is isomorphic to a subgroup of the affine group $\Aff(1,p)$.

We have $K/C = (H \times Z)/C = H/C \times ZC/C \cong H/C \times Z$. This 
implies that~$H/C$ is trivial, i.e.\ $H = C \unlhd K$.
\end{prf}

\begin{lem}
\label{TrivialCLemma}
Let the assumptions and notation be as in~{\rm \ref{Case2Section}}. 
If $H_1 = 1$, then $C = 1$ or $H$ is a $p$-group.
\end{lem}
\begin{prf}
Put $L := C_G(s)$. Then~$L$ is cyclic, $G = L \ltimes K_1$ and 
$K = H \ltimes K_1$, with $|K_1| = p$
(see Corollary~\ref{StructureOfK1} and Lemma~\ref{Case2Properties}). In 
particular,~$H$ is abelian, as it is isomorphic to a subgroup of~$L$. 

Let~$\ell$ be a prime different from~$p$ and
let $S \leq H$ denote a Sylow $\ell$-subgroup. As $|G\colon L| = p$, there is
$g \in G$ such that $S^g \leq L$. As~$L$ is abelian, we have $L \leq C_G(S^g)$.
Suppose that $C_G(S^g) = G$, i.e.\ $S^g \leq Z(G)$. Then $S \leq Z(G)$ wich
implies $S = 1$, as the core of~$H$ in~$G$ is trivial. 

Now assume that $\ell \mid |H|$. By the above, we must have
$C_G(S^g) = L$. Then $C_G( S )$ is a cyclic complement of~$K_1$ in~$G$ 
containing~$H$. As $C \leq H \leq C_G(S)$ and $C \unlhd K$, it follows that 
$C \unlhd G$ und thus $C = 1$.
\end{prf}

\begin{cor}
\label{TrivialCCorollary}
Let the assumptions and notation be as in~{\rm \ref{Case2Section}}.
If $H_1 = 1$, then $C = 1$ or $H \unlhd K$.
\end{cor}
\begin{prf}
Suppose that $C \neq 1$. Then $K = HK_1$ is a $p$-group by 
Lemma~\ref{TrivialCLemma}. 
As $K/C$ is isomorphic to a subgroup of the affine group $\Aff(1,p)$ by Galois' 
theorem (see \cite[Satz~II.3.6]{HuI}), it follows that $|K/C| = p$. Now
$C \cap C^s = 1$ and $C C^s = C \times C^s \unlhd G$ by 
Lemma~\ref{NormalSubgroupsOfK}. Thus $C^s \cong (C \times C^s)/C \leq 
K/C$, and hence~$C$ has order~$p$. Therefore, $|K| = p^2$ and hence~$K$ is
abelian, proving our claim.
\end{prf}

\subsection{The case $p = 2$.} Let us assume throughout this subsection that
$p = 2$. Then $H \unlhd K$. Here, we put $D := \cap_{g \in G} K^g$, the kernel
of the action of~$G$ on the set of right cosets of~$K$. Then~$D$ is an
elementary abelian $2$-group, as~$G$ acts faithfully on the set of right
cosets of~$H$ and as~$H$ has index~$2$ in~$K$. We
write $\bar{\ } : G \rightarrow \bar{G} := G/D$ for the canonical epimorphism.
Then~$\bar{G}$ is a faithful permutation group on~$q$ letters, i.e.\ $\bar{G}$
is isomorphic to a subgroup of~$S_q$.

\begin{lem}
\label{KIntersectT}
We have $|K \cap T| = 2$, and writing $K \cap T = \{ 1, z \}$, we have
$z \in Z(K)$. In particular, $C_G( z ) \in \{ K, G \}$.
\end{lem}
\begin{prf}
The first assertion follows from $|K\colon H| = 2$, and the second from the
fact that all $K$-conjugates of~$z$ again are in $K \cap T$.
\end{prf}

\subsubsection{Case 1} Here, we consider the case that~$K$ is normal in~$G$.
Let us keep the notation of Lemma~\ref{KIntersectT} in the following.

\begin{lem}
\label{KNormalInG}
Suppose that $K \unlhd G$. Then~$G$ is abelian and $H = 1$.
\end{lem}
\begin{prf}
In this case $z^g \in K \cap T$ for all $g \in G$, and thus $\langle z \rangle
\leq Z(G)$. Now $\langle z \rangle \cap H = 1$, as the core of~$H$ in~$G$ is
trivial. It follows that $|z| = 2$ and $K = H \times \langle z \rangle$.
This implies that $H' = K' \unlhd G$, and thus $K' = 1$, i.e.\ $K$ is abelian.
In turn, $K$ is a $2$-group, as $O_{2'}( H ) = O_{2'}( K ) \unlhd G$ (recall 
that $O_{2'}( H )$ denotes the largest normal subgroup of~$H$ of odd order). 
Thus~$K$ is the unique Sylow $2$-subgroup of~$G$.

Now let $t \in T \setminus K$. Then~$t$ is not a $2$-element as otherwise
$t \in K$. Let~$r$ be an integer such that $Q := \langle t^r \rangle$ is a
Sylow $q$-subgroup of~$G$. As~$G/K$ is cyclic, we have $G' \leq K$ and we may
thus apply Lemma~\ref{LengthsOfConjugacyClasses}. This yields
$|G \colon C_G(t)| \leq 2$. Hence $|G\colon C_G(Q)| \leq 2$ and thus
$C_G(Q) \unlhd G$. Now~$Q$ is abelian and hence $Q \leq C_G(Q)$. It follows
that $Q \unlhd C_G(Q)$ and thus $Q = O_q( C_G(Q) )$. In particular,
$Q \unlhd G$, implying that $G = Q \times K$ is abelian.
\end{prf}

\subsubsection{Case 2} Here, we consider the case that~$K$ is not normal in~$G$.
Again, we use the notation of Lemma~\ref{KIntersectT}.

\begin{prp}
\label{KNotNormal}
Suppose that $N_G(K) = K$. Then there is $M \unlhd G$ such that 
$|G\colon HM| = 2$.
\end{prp}
\begin{prf}
Let~$Z$ denote the $G$-conjugacy class of~$z$ and put $T' := T \setminus 
( Z \cup \{ 1 \} )$. By Lemma~\ref{KIntersectT}, we have $C_G(z) \in 
\{ K, G \}$. Suppose first that $C_G(z) = K$. Then $|Z| = |G\colon K| = q$, 
and thus $|T'| = q - 1$. If  $C_G(z) = G$, then $|T'| = 2(q - 1)$. In 
particular, $q - 1 \mid |T'|$.
Let~$X_1, \ldots , X_m$ denote the $G$-conjugacy classes contained in~$T'$,
numbered in such a way that $|X_1| \leq \cdots \leq |X_m|$.
Thus $|X_1| \leq q - 1$, unless $m = 1$ and $C_G(z) = G$, in which case
$|X_1| = 2(q-1)$.

Let $t \in T'$. Then $t \not\in K$, as $t \not\in \{ 1, z \} = 
K \cap T$. In particular, $\bar{t} \neq 1$, since $D \leq K$. 
Let~$X$ denote the $G$-conjugacy class of~$t$.
Then $\bar{X}$ is the $\bar{G}$-conjugacy class of~$\bar{t}$ and
$|\bar{X}|$ divides $|X|$. Consider the case that $X = X_1$.
If $|X| \leq q - 1$, then $|\bar{X}| \leq q - 1$. If
$|X| = 2(q-1)$, then $|\bar{X}|$ is a proper divisor of $2(q-1)$ by
Lemma~\ref{CentralizersInSp}(d), and thus, again, $|\bar{X}| \leq q - 1$.
Lemma~\ref{CentralizersInSp}(c) implies that $\bar{G}$ has a normal Sylow 
$q$-subgroup~$Q$. Moreover, $|\bar{G}\colon \bar{K}| = |G\colon K| = q$, and 
$\bar{K} \neq 1$, as otherwise $K = D$ would be normal in~$G$.
Thus $\bar{G}$ is a Frobenius group of order $qr$ with $r \mid (q - 1)$,
again by Lemma~\ref{CentralizersInSp}(c). 

Since~$\bar{G}$ is a Frobenius group, every non-trivial conjugacy class 
of~$\bar{G}$ has length~$q$ or~$r$, and the conjugacy classes of length~$r$ lie 
in~$Q$. 
Suppose that there is some $1 \leq j \leq m$ such that $|\bar{X}_j| = q$.
Then $|X_j| = q$ as $|\bar{X}_j|$ divides $|X_j|$
and $|X_j| \leq 2(q-1)$. Also,~$X_j$ is the unique conjugacy class of length~$q$ 
contained in $T'$. If $1 \leq i \neq j \leq m$, then $|X_i| \leq q - 1$, and
thus $|\bar{X}_i| = r$. In particular, $r \mid |X_i|$. Now $q - 1$ 
divides~$|T'|$, as we have already observed above. It follows that~$r$ 
divides~$q$, a contradiction. This shows that $\bar{T}' \subseteq Q$.

We have $z \in Z(K)$ and $D \leq K$, and thus $Z \subseteq C_G(D)$ as 
$D \unlhd G$. Moreover,~$D$ is abelian and hence $\langle Z, D \rangle 
\leq C_G(D)$. 
Now $\bar{G} = \langle \bar{T} \rangle = 
\langle \bar{Z} \cup \bar{T}' \rangle = \langle \bar{Z} \rangle$, as 
$\bar{T}' \subseteq Q$, the Frattini subgroup of~$\bar{G}$.
Thus $G = \langle Z, D \rangle \leq Z_G(D)$, i.e.\ $D \leq Z(G)$.
This implies that $D \cap H = 1$, as the core of~$H$ in~$G$ is trivial.
Hence $|H||D| = |HD| \leq |K| = 2|H|$, and so $|D| \in \{ 1, 2 \}$. 

Let~$M$ denote the inverse image of~$Q$ in~$G$. If
$D = 1$, then $|M| = q$ an thus $G = K \ltimes M$ and $|G\colon HM| = 2$
as claimed. Now suppose that $D = \langle d \rangle$ with $|d| = 2$. Then
$HM \neq G$ as otherwise $d \in HM$ and, in turn, $d \in H$. Now $G = KM$
as $\bar{G} = \bar{K}\bar{M}$, and thus $G = KM = HM \cup HMd$, i.e.\
$|G\colon HM| = 2$, as claimed.
\end{prf}

\subsection{The main result}
We can now summarize our results for envelopes of $\rcc$ loop folders of 
orders $2p$ for odd primes~$p$.

\begin{thm}
\label{MainResult}
Let $(G,H,T)$ be the envelope of an $\rcc$ loop of order $2p$, where~$p$ is an
odd prime. Then there is a subgroup $K \leq G$ with $H \leq K$ and $|G\colon K| 
= 2$ and one of the following occurs.

{\rm (a)} The group~$G$ is isomorphic to the wreath product $C_p \wr C_2$.

{\rm (b)} The group~$G$ is isomorphic to a subgroup of the affine group
$\Aff(1,p)$.

{\rm (c)} We have $G = K \times \langle a \rangle$, and~$K$ has odd order and 
is isomorphic to a subgroup of the affine group $\Aff(1,p)$.

In Cases~{\rm (b)} and~{\rm (c)}, $\langle T \cap K \rangle$ is a normal 
subgroup of~$G$ of order~$p$.
The Cases~{\rm (a)},~{\rm (b)} and~{\rm (c)} are disjoint.
\end{thm}
\begin{prf}
The first statement follows from Lemmas~\ref{KNormalInG} and 
Proposition~\ref{KNotNormal} (with~$q$ replaced by~$p$). In particular, we are
in the situation of Subsection~\ref{CaseQEqal2}.

In the following, we resume to the notation introduced at the beginning of 
Subsection~\ref{Generalities}.
Suppose that~$G$ is not isomorphic to the wreath product $C_p \wr C_2$. By
Proposition~\ref{Case1}, we may assume that we are in the situation 
of~\ref{Case2Section}. Corollary~\ref{HNormalInK} implies that~$H$ is not a 
normal subgroup of~$K$. Hence $H_1 = 1$ and $C = 1$ by 
Corollaries~\ref{Case2Special} and~\ref{TrivialCCorollary}. In particular,
$|K_1| = p$ by Corollary~\ref{StructureOfK1}.

If $C_G(K_1) = K_1$, then $G/K_1$ injects into the automorphism group of~$K_1$,
and thus~$G$ is as in~(b). Assume now that $K_1 \lneq C_G(K_1)$. As $K = HK_1$
by Lemma~\ref{FirstProperties}, we have $C_H(K_1) \leq C$,
and thus $H \cap C_G(K_1) = C_H(K_1) = 1$. Hence
$|HC_G(K_1)| = |H||C_G(K_1)| > |H||K_1| = |K|$, and thus $HC_G(K_1) = G$ and
$|C_G(K_1)| = 2|K_1|$. It follows that $C_G(K_1) = K_1 \times \langle a \rangle$
for some $a \in G \setminus K$ of order~$2$. As $K \unlhd G$, and 
$\langle a \rangle = O_2(C_G(K_1)) \unlhd G$, we have
$G = K \times \langle a \rangle$. Now~$K$ acts faithfully on the set of 
$H$-cosets in~$K$, and thus $K = HK_1$ is isomorphic to a subgroup of the 
\enlargethispage{2em}
affine group $\Aff(1,p)$. Finally,~$K$ has odd order since $G/K_1$ is cyclic by 
Lemma~\ref{Case2Properties}.
\end{prf}

\section{The rcc loops of twice prime order}\label{6}
\markright{RCC LOOPS OF TWICE PRIME ORDER}

Let~$p$ be an odd prime. In this section we determine the number of isomorphism
classes of \rcc loops of order~$2p$. Let~$\mathcal{L}$ denote such a loop and 
let $(G,H,T)$ be its envelope. By numbering the elements of~$\mathcal{L}$ by
the integers $1, \ldots , 2p$, where~$1$ numbers the identity element 
of~$\mathcal{L}$, we may and will view~$G$ as a subgroup 
of~$S_{2p}$, and~$H$ as the stabilizer in~$G$ of~$1$. If~$\mathcal{L}_1$
and $(G_1,H_1,T_1)$ is another such configuration, then~$\mathcal{L}$ and
$\mathcal{L}_1$ are isomorphic as loops, if and only if there is an element 
of~$S_{2p}$, conjugating $(G,H,T)$ to $(G_1,H_1,T_1)$. The isomorphism types
of the right multiplication groups arising in \rcc loops of order $2p$ have
been described in Theorem~\ref{MainResult}. For each of these groups~$G$ we
have to determine their embeddings into $S_{2p}$ up to conjugation. This
will yield the possible pairs $(G,H)$ to be considered. For each of these pairs
we have to determine the normalizer~$N$ in $S_{2p}$ of~$G$ and~$H$, and then
find the distinct $N$-orbits of $G$-invariant transversals~$T$ for $H\backslash G$
such that $1 \in T$ and $\langle T \rangle = G$. We will refer to the three 
different
types of~$G$ in Theorem~\ref{MainResult}(a),~(b), and~(c) as Case~(a),~(b) 
and~(c), respectively.

We begin with some preliminary results. As usual, the largest normal 
$p$-subgroup of a finite group~$U$ is denoted by~$O_p(U)$, and its largest
normal subgroup of odd order by $O_{2'}(U)$.

\begin{lem}
\label{NormalizersInS2p}
Let
$\pi_1, \alpha \in S_{2p}$ be defined by
$\pi_1 := (1,2, \ldots ,p)$
and
$\alpha := (1,p+1)(2,p+2) \cdots (p,2p)$.
Put $\pi_2 := \pi_1^\alpha = (p+1, p + 2, \ldots , 2p)$.
Let $\nu_1 \in S_p$ be an element of order $p - 1$ such that
$N_{S_p}( \langle \pi_1 \rangle ) = \langle \pi_1, \nu_1 \rangle$.
Put $\nu_2 := \nu_1^\alpha$ and $\nu := \nu_1 \nu_2$.

{\rm (a)} Let $G := \langle \pi_2, \alpha \rangle$. Then 
$Z(G) = \langle \pi_1\pi_2 \rangle$, $G' = \langle \pi_1^{-1}\pi_2 \rangle$
and $G = C_{S_{2p}}(Z( G )) \cong C_p \wr C_2$. Put $N := N_{S_{2p}}( G )$.
Then $N = N_{S_{2p}}(Z( G )) = \langle \nu \rangle \ltimes G$. 

{\rm (b)} Let $U \leq N$ with $O_p(U) = Z(G)$ and $U \not\leq G$. Then there is
$n \in N$ such that $N_{S_{2p}}( U^n ) = A \times \langle \alpha \rangle$ with 
$A = \langle \nu \rangle \ltimes Z(G)$.
\end{lem}
\begin{prf}
(a) The statements about~$G$ are trivially verified. From 
$G = C_{S_{2p}}(Z( G ))$ we conclude that $G \unlhd N_{S_{2p}}(Z( G ))$, and
thus $N = N_{S_{2p}}(Z( G ))$. Moreover, $N/G$ is isomorphic to a subgroup of 
$\Aut( Z( G ) )$, which is a cyclic group of order $p - 1$. Now 
$\langle \nu \rangle \cap G = 1$, as the elements in~$G$ fixing the set 
$\{ 1, \ldots , p \}$ have order divisible by~$p$. Also,~$\nu$ 
normalizes~$G$, and thus $N = \langle \nu \rangle \ltimes G$. 

(b) From $O_p(U) = Z(G)$ we conclude that $N_{S_{2p}}( U ) \leq 
N_{S_{2p}}(Z( G )) = N$, and thus $N_{S_{2p}}( U ) = N_N( U )$. Let~$V$ denote
a complement to~$Z(G)$ in~$U$, and let~$W$ be a Hall $p'$-group of~$N$ 
containing~$V$ (see \cite[Hauptsatz~VI.1.7]{HuI}). Then~$W$ is a complement to 
$O_p( N )$ in~$N$. As~$\alpha$ 
centralizes~$\nu$, we have $\langle \nu, \alpha \rangle = 
\langle \nu \rangle \times \langle \alpha \rangle$, and thus 
$\langle \nu \rangle \times \langle \alpha \rangle$ is another complement to 
$O_p( N )$ in~$N$. As all such complements are conjugate in~$N$, there is 
$n \in N$ such that $W^n = \langle \nu \rangle \times \langle \alpha \rangle$
and $U^n = Z(G)V^n$. By replacing~$U$ with $U^n$, we may assume that 
$V \leq W = \langle \nu \rangle \times \langle \alpha \rangle$.
In particular,~$W$ is abelian. It follows that~$WZ(G)$ normalizes~$U = VZ(G)$. 
As $U \not\leq G$, there is an element $\nu^i\alpha^j \in V$ such that 
$\nu^i \neq 1$. Then $[\nu^i\alpha^j,\pi_1] \not\in Z(G)$.
In particular,~$U$ is not normal in~$N$. As~$WZ(G)$ has index~$p$ in~$N$, 
we conclude that $N_N(U) = WZ(G)$, which proves our claim.
\end{prf}

Let $n, d$ be positive integers, and let $\zeta \in S_n$ denote an 
$n$-cycle. Let us put
\begin{equation}
\label{Involutions1}
I_{n,d} := 
|\{ \tau \in S_n \mid \tau^2 = 1, \tau \zeta^d = \zeta^d \tau \}|
\end{equation}
and
\begin{equation}
\label{Involutions2}
I_n := I_{n,n}.
\end{equation}
Thus~$I_{n,d}$ is one more than the number of involutions in 
$C_{S_n}( \zeta^d )$ and~ $I_n$ is one more than the number of involutions 
in~$S_n$. Notice that the definition of~$I_{n,d}$ does not depend on the chosen 
$n$-cycle~$\zeta$, as all $n$-cycles are conjugate in~$S_n$.

It is not difficult to derive a formula for~$I_{n,d}$, where the formula 
for~$I_n$ is certainly well known. In the following result, 
$n \text{\rm\ mod\ } 2 \in \{ 0, 1 \}$ denotes the remainder of the division 
of~$n$ by~$2$.
\begin{lem}
\label{IndFormula}
Let $n, d, e$ and~$f$ be positive integers such that $d \mid n$ and 
$\mbox{\rm gcd}(e,n/d) = 1$. Then $I_{n,de} = I_{n,d}$ and $I_{n,f} = 
I_{n,\gcd( n, f )}$.
Moreover, we have
$$I_{n,d} = \sum_{k = 0}^{\lfloor d/2 \rfloor} 
\frac{d! (n/d)^k (2 - (n/d \text{\rm\ mod\ } 2))^{d - 2k}}{2^kk!(d - 2k)!}.$$
In particular,
$$I_n = \sum_{k = 0}^{\lfloor n/2 \rfloor}
\frac{n!}{2^kk!(n - 2k)!}.$$
\end{lem}
\begin{prf}
Let $\zeta \in S_n$ be an $n$-cycle. As~$e$ is relatively prime to~$n/d$, we
have that $\zeta^{de} = (\zeta^d)^e$ is the 
product of~$d$ cycles of length~$n/d$. In particular,~$\zeta^d$ and $\zeta^{de}$
are conjugate in~$S_n$ and thus $I_{n,de} = I_{n,d}$. Writing $f = de$ with
$d = \gcd( n, f)$ and $e = f/\gcd( n, f)$, we obtain $I_{n,f} = 
I_{n,\gcd( n, f )}$, as $f/\gcd( n, f)$ and $n/\gcd( n, f)$ are relatively prime.

By definition, $I_{n,d}$ equals the number of elements 
$\tau \in C_{S_n}( \zeta^d )$ with $\tau^2 = 1$.
The structure of $C_{S_n}( \zeta^d )$ is well known; it is a wreath product
isomorphic to $C_{n/d} \wr S_d$, where~$C_{n/d}$ denotes a cyclic group of
order~$n/d$. We view the elements of $C_{S_n}( \zeta^d )$ as $(d + 1)$-tuples
$(\mu;c_1, \ldots , c_d )$, where each~$c_i$ lies in one of the~$d$ cycles 
of~$\zeta^d$, and where~$\mu \in S_d$ permutes the numbers 
$\{ 1, \ldots , d \}$. We have 
$$(\mu; c_1, \ldots , c_d )^2 = (\mu^2; c_1c_{1\mu^{-1}}, c_2c_{2\mu^{-1}}, 
\ldots , c_dc_{d\mu^{-1}}).$$
Let $\tau := (\mu; c_1, \ldots , c_d ) \in C_{S_n}( \zeta^d )$
satisfy $\tau^2 = 1$. Then $\mu^2 = 1$ and $c_{i\mu} = c_i^{-1}$ for all 
$1 \leq i \leq d$. Suppose that~$\mu$ is a product of exactly~$k$ 
transpositions for some $0 \leq k \leq \lfloor d/2 \rfloor$. Then $c_j = 
c_i^{-1}$, if $(i,j)$ is a transposition of~$\mu$, and $c_i^2 = 1$ if~$i$
is a fixed point of~$\mu$. This way, a fixed~$\mu$ gives rise to $(n/d)^k
(2 - (n/d \text{\rm\ mod\ } 2))^{d - 2k}$ elements $\tau \in C_{S_n}( \zeta^d )$
with $\tau^2 = 1$. The centraliser of~$\mu$ in~$S_d$ has order $2^k k! (d - 2k)!$,
yielding our formula for $I_{n,d}$. The one for~$I_n$ follows from this by 
putting $d = n$.
\end{prf}

\begin{prp}
\label{Casea}
There are exactly 
$$I_{p-1} - 1 + \frac{1}{p-1}\sum_{d=1}^{p-1} I_{p-1,d}$$ 
distinct isomorphism types of $\rcc$ loops with
multiplication group~$G$ as in {\rm Case~(a)}. 
\end{prp}
\begin{prf}
Let $(G,H,T)$ denote the envelope of an $\rcc$ loop of order~$2p$ with~$G$ as in
Case~(a), i.e.\ $G$ is isomorphic to the wreath product $C_p \wr C_2$. In this 
case,~$H$ is cyclic of order~$p$. By numbering the right cosets of~$H$ in~$G$ 
from~$1$ to~$2p$, we obtain an embedding $G \rightarrow S_{2p}$, and we 
identify~$G$ with its image in~$S_{2p}$ from now on. Let~$\pi_1$,~$\alpha$ 
and~$\pi_2$ be defined as in Lemma~\ref{NormalizersInS2p}.
We may choose the numbering of the cosets of~$H$ in~$G$ in such a way that
$H = \langle \pi_2 \rangle$ and $G = \langle \pi_2, \alpha \rangle$. From 
Lemma~\ref{NormalizersInS2p}(a) we obtain $Z(G) = \langle \pi_1\pi_2 \rangle$, 
$G' = \langle \pi_1^{-1}\pi_2 \rangle$ and $G = C_{S_{2p}}(Z( G ))$. Also, 
$N := N_{S_{2p}}( G )$ equals $\langle \nu \rangle \ltimes G$
with~$\nu$ as in Lemma~\ref{NormalizersInS2p}.
Observe that~$N$ normalises~$H$.

Let~$\mathcal{T}$ denote the set of $G$-invariant transversals 
for~$H\backslash G$ containing~$1$. Put $K := \langle \pi_1, \pi_2 \rangle = 
O_p(G)$ and let ~$\mathcal{T}_1$ denote the set of $G$-invariant transversals 
for~$H\backslash K$ containing~$1$. Let
$t \in G \setminus K$. Then $|C_G( t )| = 2p$ and thus~$t$ lies in a conjugacy
class of length~$p$. As every conjugacy class of~$G$ lies in some coset of~$G'$,
we find that $G't$ is the conjugacy class of~$G$ containing~$t$. Hence if $T \in
\mathcal{T}$, 
we have $T = (K \cap T) \cup G't$ for some $t \in G \setminus K$, and $K \cap T
\in \mathcal{T}_1$.
Conversely, if $T_1 \in \mathcal{T}_1$, 
and if~$t$ is any element of $G \setminus K$,
then $T_1 \cup G't \in \mathcal{T}$.

As $K = H \times H^\alpha$, we have $K = \cup_{0 \leq j \leq p - 1} H \pi_1^j$.
A transversal for $H\backslash K$ contains exactly one element of each coset
$H \pi_1^j$, $0 \leq j \leq p - 1$. As we insist that our transversals contain
the trivial element, a transversal~$T_1$ for $H\backslash K$ determines a map
$\tau : \{ 1, \ldots , p - 1 \} \rightarrow \{ 0, 1, \ldots , p - 1 \}$ such
that 
\begin{equation}
\label{T1}
T_1 = 
\{ \pi_2^{j \tau} \pi_1^j \mid 1 \leq j \leq p - 1 \} \cup \{ 1 \}. 
\end{equation}
Conjugating the element
$\pi_2^{j \tau} \pi_1^j \in T_1 \setminus \{ 1 \}$ by~$\alpha$, we obtain
$\pi_1^{j \tau} \pi_2^j = \pi_2^j \pi_1^{j \tau}$. If~$T_1$ is $G$-invariant,
we must have, firstly, that $j \tau \neq 0$ and, secondly, that 
$\pi_2^j \pi_1^{j \tau} \in T_1 \setminus \{ 1 \}$ for all $1 \leq j \leq p - 1$.
The latter condition implies that $\pi_2^{ j \tau^2} \pi_1^{j \tau} = 
\pi_2^j \pi_1^{j \tau}$ for all $1 \leq j \leq p - 1$, and thus $\tau^2 = 1$.
In particular,~$\tau$ is a permutation of order at most~$2$ of the set 
$\{ 1, \ldots , p - 1 \}$. Conversely, if~$\tau$ is a permutation of the latter
set with $\tau^2 = 1$, then $T_1$ defined by~(\ref{T1}) lies in~$\mathcal{T}_1$.
In particular, $|\mathcal{T}_1| = I_{p-1}$.
As the 
number of conjugacy classes of~$G$ in $G \setminus K$ equals~$p$, we conclude
from 
$$\mathcal{T} = \{ T_1 \cup G't \mid T_1 \in \mathcal{T}_1, t \in G\setminus K \},$$
that
$$|\mathcal{T}| = pI_{p-1}.$$

We next determine the number of $N$-orbits on~$\mathcal{T}$. This is the
same as the number of $\langle \nu \rangle$-orbits on~$\mathcal{T}$. To compute
this number, put
$$\mathcal{T}' := \{ T_1 \cup G'\alpha \mid T_1 \in \mathcal{T}_1 \} \subseteq 
\mathcal{T}.$$
Observe that $\mathcal{T}_1$ is $\langle \nu \rangle$-invariant, as~$\nu$ 
normalises~$H$. In addition,~$\nu$ centralises~$\alpha$, and thus
$\mathcal{T}'$ is $\langle \nu \rangle$-invariant as well.
As~$Z(G)$ is a set of representatives for the set of right cosets of~$G'$ 
in~$K$, every conjugacy class of~$G$ contained in $G \setminus K$ is of the form
$G'z\alpha$ for some $z \in Z(G)$. As $\langle \nu \rangle$ acts transitively on
$Z(G) \setminus \{ 1 \}$, we conclude that
every orbit of~$\langle \nu \rangle$ on $\mathcal{T} \setminus 
\mathcal{T}'$ has length~$p - 1$, and thus there are exactly $I_{p-1}$ such 
orbits. We are thus left with the determination of the number of 
$\langle \nu \rangle$-orbits on~$\mathcal{T}'$, which is the same as the number 
of $\langle \nu \rangle$-orbits on~$\mathcal{T}_1$. By the 
Burnside-Cauchy-Frobenius lemma, the latter number equals
$$\frac{1}{p-1}\sum_{d = 1}^{p-1} \chi_d,$$
where~$\chi_d$ is the number of fixed points of~$\nu^d$ on~$\mathcal{T}_1$.
The action of~$\langle \nu \rangle$ on~$K$ determines a $(p - 1)$-cycle~$\zeta$
on the set $\{ 1, \ldots , p - 1 \}$ such that $\nu^{-1} \pi_i^j \nu = 
\pi_i^{j\zeta}$ for $i = 1, 2$ and all $1 \leq j \leq p - 1$. Now let
$T_1 \in \mathcal{T}_1$ be given by~(\ref{T1}) with respect to 
$\tau \in S_{p-1}$ with $\tau^2 = 1$. Then~$T_1$ is fixed by~$\nu^d$, if and
only if $\zeta^d$ centralises~$\tau$. Thus $\chi_d = I_{p-1,d}$. 

It remains to determine those $N$-orbits on~$\mathcal{T}$ containing 
transversals that generate~$G$. Let $T \in \mathcal{T}$ such that 
$\langle T \rangle \neq G$. Then $\langle T \rangle$ is a normal subgroup of~$G$
of index~$p$. Thus $G' \leq \langle T \rangle$ and $T = G' \cup G'z\alpha$ for 
some $z \in Z(G)$. Since $\langle T \rangle \neq G$, we must have $z = 1$,
i.e.\ $\langle T \rangle = T = G' \cup G'\alpha$. As this is $N$-invariant,
our result follows.
\end{prf}

\begin{prp}
\label{Casebc}
Write $p - 1 = 2^nr$ with positive integers $n$ and $r$ and with~$r$ odd. Then 
there are exactly $p - r - 1$ distinct isomorphism types of $\rcc$ loops with
multiplication group~$G$ as in {\rm Case~(b)}, and there are exactly~$r$ 
isomorphism types of $\rcc$ loops with multiplication group~$G$ as in 
{\rm Case~(c)}.
\end{prp}
\begin{prf}
Let $(G,H,T)$ denote the envelope of an $\rcc$ loop of order~$2p$ with~$G$ as in
Case~(b) or~(c). If $H = 1$, then $T = G$ is a group of order~$2p$, which is
non-abelian in Case~(b), and cyclic in Case~(c). In each case, we obtain a 
unique isomorphism class of $\rcc$ loops. 

Thus let us assume that $H \neq 1$ in the following. As in the proof of 
Proposition~\ref{Casea}, we identify~$G$ with its image in~$S_{2p}$ through an
embedding obtained by numbering the right cosets of~$H$ in~$G$ from~$1$ to~$2p$. 
Put $P := O_p(G)$, the unique Sylow $p$-subgroup of~$G$.
Let $\pi_1$,~$\alpha$ and~$\pi_2$ be defined as in Lemma~\ref{NormalizersInS2p}.
We may choose the numbering of the cosets of~$H$ in~$G$ in such a way that 
$P = \langle \pi_1\pi_2 \rangle$, and that $Z(G) = \langle \alpha \rangle$ in 
Case~(c).
Put $N := N_{S_{2p}}(G)$. We now apply Lemma~\ref{NormalizersInS2p}(b) with 
our~$G$ taking the role of~$U$ of that lemma. As $H \neq 1$, we have $G \not\leq
\langle \pi_2, \alpha \rangle$, and thus, replacing~$G$ by a suitable conjugate 
within $N_{S_{2p}}( P )$, we find that $N = A \times \langle \alpha \rangle$, with 
$A \cong \Aff(1,p)$. We have $A = L \ltimes P$, with~$L$ cyclic of order $p - 1$. 

Assume that~$G$ is as in Case~(b). Then $G \cap L$ is a complement to~$P$ 
in~$G$. As all such complements are conjugate in~$G$ by Schur's theorem
(see \cite[Satz I.17.5]{HuI}), we
may assume that $H \leq L$. In particular, $G \leq A$, and~$H$ is $N$-invariant. 
Let~$T$ be a $G$-invariant transversal for $H\backslash G$. Then 
$P \subseteq T$ by Theorem~\ref{MainResult}. Let $\tau \in T \setminus P$. Then 
$|C_G(\tau)| = 2|H|$ and thus $T \setminus P$ consists of the $G$-conjugacy 
class containing~$\tau$. If, moreover, $G = \langle P, \tau \rangle$, we have 
$2p|H| = |G| = p|\tau|$ and~$\tau$ has even order larger than~$2$. Every 
element~$\tau'$ which is conjugate to~$\tau$ in~$A$ gives rise to an isomorphic 
loop with multiplication group $\langle P, \tau' \rangle$, as in Case~(b).
It follows that the isomorphism types of $\rcc$ loops with a multiplication
group as in Case~(b) equals the number of $A$-conjugacy
classes of elements of~$A$ of even order larger than~$2$. As~$A$ has
$(p - r - 2)p$ such elements, the result follows.

Assuming now that~$G$ is as in Case~(c), we have 
$G = K \times \langle \alpha \rangle$, with $K = O_{2'}(G)$, and thus 
$K \unlhd N$. In turn, $K \leq A$ as every Sylow subgroup of~$K$ is conjugate to 
a subgroup of~$A$. 
Again,~$H$ is $N$-invariant. Let~$T$ be a $G$-invariant transversal for 
$H\backslash G$. As in Case~(b), we have $T = P \cup C$, where~$C$ is a 
$G$-conjugacy class of an element $\tau \in G \setminus K$. Every
element~$\tau'$ in the $A$-conjugacy class containing~$\tau$ gives rise to an
isomorphic loop with multiplication group $\langle P, \tau' \rangle$.
Now $\tau = \tau_1 \alpha$ for some~$\tau_1 \in K$. It follows that the 
isomorphism types of $\rcc$ loops 
with a multiplication group as in Case~(c) equals the number of $A$-conjugacy
classes of elements of~$A$ of odd order different from~$p$. All these
elements lie in the unique subgroup of~$A$ of order $pr$, and thus there
are $(r-1)p$ non-trivial such elements. As the trivial element yields a group,
the result follows.
\end{prf}

We summarise our results in the following theorem.
\begin{thm}
\label{NumberOfAllLoops}
Let~$p$ be a prime. Then the number of isomorphism types of $\rcc$ loops of
order $2p$ (including groups) equals
\begin{equation}
\label{FormulaForNumberOfAllLoops}
p - 2 + I_{p-1} + \frac{1}{p-1}\sum_{d=1}^{p-1} I_{p-1,d}.
\end{equation}
\end{thm}
\begin{prf}
Every loop of order~$4$ is a group. As $I_{1,1} = I_1 = 1$, 
formula~(\ref{FormulaForNumberOfAllLoops}) holds for $p = 2$. For odd~$p$
it follows from Propositions~\ref{Casea} and~\ref{Casebc}, as the cases in
Theorem~\ref{MainResult} are disjoint.
\end{prf}

The table below contains the numbers obtained by evaluating
formula~(\ref{FormulaForNumberOfAllLoops}) for small values of~$p$. These 
numbers have also been obtained for $p \leq 13$ in the PhD-thesis of the first 
author~\cite{Artic} by different methods.

$$\begin{array}{c|cccccccc}
 p & 2 & 3 & 5 & 7 & 11 & 13 & 17 & 19 \\ \hline
{\rm (\ref{FormulaForNumberOfAllLoops})} & 2 & 5 & 18 & 99 & 10\,489
& 151\,973 & 49\,096\,721 & 1\,052\,729\,657
\end{array}
$$

\medskip

\noindent
One of the referees has kindly pointed out that 
formula~(\ref{FormulaForNumberOfAllLoops}) evaluates to an integer, even if~$p$ 
is not a prime (and larger than~$1$). This follows from the fact that for general positive integers
$n, d$, the number $I_{n,d}$ equals the number of fixed points of the element 
$\zeta^d$ on the set $\{ \tau \in S_n \mid \tau^2 = 1 \}$, where the $n$-cycle
$\zeta$ acts by conjugation. Thus, by the Burnside-Cauchy-Frobenius lemma, the 
number of orbits of $\langle \zeta \rangle$ on $\{ \tau \in S_n \mid \tau^2 = 1 \}$
equals $1/n \sum_{d = 1}^{n} I_{n,d}$, so that this number is an integer.

\section{A series of examples}\label{7}
\markright{A SERIES OF EXAMPLES}

According to Theorem~\ref{MainResult}, the right multiplication group of an 
$\rcc$ loop of order $2p$, where~$p$ is an odd prime, is soluble. This is no
longer the case for right multiplication groups of $\rcc$ loops of order $pq$, 
where~$p$ and~$q$ are distinct primes. An example is given 
in~\cite[Table~B.7]{Artic} of an $\rcc$ loop of order~$15$ with right 
multiplication group isomorphic to~$\GL(2,4)$. This fits into an infinite
series of examples.

\begin{prp}
\label{GL2Examples}
Let~$q$ be a power of a prime with $q > 2$. Then there is an $\rcc$ loop of 
of order $q^2 - 1$ and right multiplication group isomorphic to $\GL(2,q)$.
\end{prp}
\begin{prf}
Let $G := \GL(2,q)$, acting from the right on $\mathbb{F}_q^{1 \times 2}$, and 
let
$$H := \left\{ 
\left( \begin{array}{cc} \alpha & 0 \\ \beta & 1 \end{array} \right)
\mid \alpha \in \mathbb{F}_q^*, \beta \in \mathbb{F}_q \right\}.$$
Let $Z := Z(G)$ denote the set of scalar matrices in~$G$ and let~$C$ be
a $G$-conjugacy class of elements of order $q^2 - 1$, i.e.\ the elements of~$C$
are Singer cycles. Then $|C_G(t)| = q^2 - 1$ for all $t \in C$;
in particular $|C| = q(q - 1)$. Now put
$$T := C \cup Z.$$
We claim that~$T$ is a $G$-invariant transversal for $H\backslash G$. Clearly,
$T$ is $G$-invariant and $|T| = q^2 - 1 = |G\colon H|$. Let $t, t' \in C$.
We have to show that $t't^{-1} \in H$ if and only if $t = t'$. To see this, 
first observe that 
$|C_G(t)||H| = |G|$ and that $C_G(t) \cap H = 1$, as $|C_G(t) \cap H|$ divides 
$\mbox{\rm gcd}( |C_G(t)|, |H| ) = q - 1$, and the only elements in~$C_G(t)$ 
of order dividing $q - 1$ are the elements of~$Z$. We conclude that 
$G = C_G(t)H$. It follows that there is $h \in H$ with $t' = h^{-1}th$. Put 
$h' := t't^{-1} = h^{-1}tht^{-1}$. Thus 
\begin{equation}
\label{i}
t^{-1}hh' = ht^{-1}.
\end{equation}
Now assume that $h' \in H$.
As $\mbox{\rm det}(h') = \mbox{\rm det}(t't^{-1}) = 1$, we have
$$h' = \left( \begin{array}{cc} 1 & 0 \\ \gamma & 1 \end{array} \right)$$
for some $\gamma \in \mathbb{F}_q$. Let
$$h = \left( \begin{array}{cc} \alpha & 0 \\ \beta & 1 \end{array} \right)$$
with 
$\alpha \in \mathbb{F}_q^*$ and $\beta \in \mathbb{F}_q$, and let
$$t^{-1} = \left( \begin{array}{cc} a & b \\ c & d \end{array} \right)$$
with $a, b, c, d \in \mathbb{F}_q$. Then
$$t^{-1}hh' = \left( \begin{array}{cc} * & b \\ * & d \end{array} \right),$$
and 
$$ht^{-1} = \left( \begin{array}{cc} * & \alpha b \\ * & \beta b + d \end{array} \right),$$
where we do not need to specify the entries in the first columns of $t^{-1}hh'$ 
respectively $ht^{-1}$.  As~$t$ acts irreducibly on the natural vector space 
$\mathbb{F}_q^{1 \times 2}$ for~$G$, we conclude that $b \neq 0$. 
Equation~(\ref{i}) yields  $\alpha = 1$ and $\beta = 0$, i.e.\ $h = 1$, and 
thus $t = t'$.
If $z, z' \in Z$, then $z' z^{-1} \in H$ if and only if $z = z'$. Now
let $z \in Z$ and $t \in C$ and assume that $tz^{-1} \in H$. Then $t \in HZ$;
but $|HZ| = q(q-1)^2$, whereas $|t| = q^2 - 1 \nmid q(q-1)^2$, a contradiction.

Finally, it is easy to check that $\langle T \rangle = G$, by a direct 
computation if $q = 3$, and by using the fact that $G/Z$ is almost simple if 
$q \neq 3$. This completes the proof.
\end{prf}

\section*{Acknowledgements}

We thank Alice Niemeyer for her support and her interest in this work.
We also thank Barbara Baumeister for introducing us to the fascinating topic
of invariant transversals. Finally, we are very much indebted to the anonymous 
referees for several suggestions which improved the exposition of this paper, 
and also for drawing our attention to related work.

\end{document}